\documentclass[12pt]{amsart}
\usepackage{amscd,amsfonts,amsmath,amsxtra,amssymb}
\usepackage[all]{xy}
\usepackage{hyperref}
\usepackage[french]{babel}
\usepackage[T1]{fontenc}

\newtheorem{sub}{}[section]
\newtheorem{subsub}{}[sub]
\font\tte=cmbsy10

\def\ov#1{\overline{#1}}
\def\codim{\mathop{\rm codim}\nolimits}
\def\coker{\mathop{\rm coker}\nolimits}
\def\Hom{\mathop{\rm Hom}\nolimits}
\def\Ext{\mathop{\rm Ext}\nolimits}
\def\Aut{\mathop{\rm Aut}\nolimits}
\def\imm{\mathop{\rm Im}\nolimits}

\def\lra{\longrightarrow}
\def\som{\mathop{\hbox{$\displaystyle\bigoplus$}}\limits}
\def\sigg{\mathop{\hbox{$\displaystyle\sum$}}\limits}

\def\supp{\mathop{\hbox{$\sup$}}\limits}

\def\timex{\setbox250=\hbox{$\times$}\hskip 5pt\Box\hskip -0.75em
{\raise 1.5pt\vbox{\box250}}\hskip 5pt}
\def\paragra{{\tte \char120}}
\def\para{\paragra~\hskip -2pt}
\def\subsubsubsection#1#2{\noindent{\bf \ #1} \ \ {\it #2}}
\def\hfl#1#2{\smash{\mathop{\hbox to 12mm{\rightarrowfill}}
\limits^{\scriptstyle#1}_{\scriptstyle#2}}}

\def\pline#1{<\hskip-3.5pt#1\hskip-3.5pt>}

\def\m#1{{\hbox{$#1$}}}

\def\fexcs{fibr\'es exceptionnels }
\def\ot{\otimes}

\newcommand{\C}{{\mathbb C}}

\renewcommand{\P}{{\mathbb P}}

\newcommand{\ke}{{\mathcal E}}
\newcommand{\kf}{{\mathcal F}}
\newcommand{\kg}{{\mathcal G}}
\newcommand{\kh}{{\mathcal H}}

\newcommand{\kk}{{\mathcal K}}

\newcommand{\ko}{{\mathcal O}}

\newcommand{\ku}{{\mathcal U}}
\newcommand{\kv}{{\mathcal V}}
\newcommand{\kw}{{\mathcal W}}

\def\og{\leavevmode\raise.3ex\hbox{$\scriptscriptstyle\langle\!\langle$}}
\def\fg{\leavevmode\raise.3ex\hbox{$\scriptscriptstyle\,\rangle\!\rangle$}}

\def\nsp{\lbrace 0\rbrace}
\def\eva#1#2{#1\ot\Hom(#1,#2)\lra #2}
\def\coeva#1#2{#1\lra #2\ot\Hom(#1,#2)^*}
\begin{document}
\newtheorem{xprop}{Proposition}[section]
\newtheorem{xlemm}[xprop]{Lemme}
\newtheorem{xtheo}[xprop]{Th\'eor\`eme}
\newtheorem{xcoro}[xprop]{Corollaire}
\newtheorem{quest}{Question}
\newtheorem{defin}{D\'efinition}

\def\refname{R\'ef\'erences}
\def\contentsname{Sommaire}
\def\proofname{D\'emonstration}
\def\abstractname{R\'esum\'e}

\markboth{\small Jean-Marc Dr\'ezet}{\small Quotients alg\'ebriques par des 
groupes non r\'eductifs}

\title[{\tiny Vari\'et\'es de modules de complexes}]
{Quotients alg\'ebriques par des groupes non r\'eductifs\break
et vari\'et\'es de modules de complexes}
\author[{J.M. Dr\'{e}zet}]{Jean-Marc Dr\'ezet}
\address{
Institut de Math\'ematiques de Jussieu,
Case 247,
4 place Jussieu,
F-75252 Paris, France}
\email{drezet@math.jussieu.fr}
\urladdr{http://www.math.jussieu.fr/$\sim$drezet}

\maketitle
\tableofcontents

\section{Introduction}

\begin{sub}{\bf Le probl\`eme des quotients alg\'ebriques par des groupes non
r\'eductifs}\end{sub}

Soit $G$ un groupe alg\'ebrique lin\'eaire de radical unipotent $H$. On
suppose qu'il existe un sous-groupe alg\'ebrique r\'eductif \m{G_{red}} de $G$
dont l'inclusion dans $G$ induit un isomorphisme \ \m{G_{red}\simeq G/H}.
Soit $Y$ une vari\'et\'e alg\'ebrique projective munie d'une action
alg\'ebrique de $G$, et $L$ un $G$-fibr\'e en droites tr\`es ample sur $Y$.
On dit qu'un point $y$ de $Y$ est {\em semi-stable} relativement \`a $L$
s'il existe un entier \ \m{k>0} \ et une section $G$-invariante $s$ de \m{L^k}
telle que \ \m{s(y)\not =0}. Soit \m{Y^{ss}(L)} l'ouvert $G$-invariant de Y
constitu\'e des points $G$-semi-stables relativement \`a $L$. La construction 
d'un {\em bon quotient} \ \m{Y^{ss}(L)//G} \ est possible dans le cas o\`u
\ \m{H=\nsp}, le groupe $G$ \'etant dans ce cas r\'eductif (cf \cite{mumf},
\cite{news}). Le cas o\`u $G$ n'est pas r\'eductif est plus difficile, et
a \'et\'e abord\'e par A. Fauntleroy dans \cite{faunt}. On doit consid\'erer
un ouvert $G$-invariant plus petit de \m{Y^{ss}(L)} (en imposant des
conditions suppl\'ementaires qui d\'ependent essentiellement de l'action de
$H$) et les quotients obtenus sont en g\'en\'eral seulement des {\em quotients
cat\'egoriques}. De plus, la d\'efinition de l'ouvert \`a quotienter est peu
explicite. Ces restrictions s'expliquent sans doute par la grande 
g\'en\'eralit\'e des probl\`emes trait\'es dans \cite{faunt}. On propose ici
une d\'efinition l\'eg\`erement diff\'erente de la semi-stabilit\'e :

\bigskip

\begin{defin}
On dit qu'un point $y$ de $Y$ est {\em $G$-semi-stable} (resp. 
{\em $G$-stable})
relativement \`a $L$ si tout point de l'orbite \m{Hy} est 
\m{G_{red}}-semi-stable
(resp. \m{G_{red}}-stable) relativement \`a $L$ (vu comme un 
\m{G_{red}}-fibr\'e en droites).
\end{defin}

\bigskip

Il est clair que les points semi-stables relativement \`a $L$ le sont aussi
au sens de la d\'efinition 1. Cette d\'efinition me semble plus explicite,
car les points $G_{red}$-semi-stables peuvent en g\'en\'eral \^etre
d\'etermin\'es \`a l'aide de crit\`eres num\'eriques (cf. \cite{mumf}).

\vskip 1.5cm

\begin{sub}{\bf Espaces de complexes}\end{sub}

On s'int\'eresse dans cet article \`a un type particulier d'action. Soient $X$
une vari\'et\'e alg\'ebrique projective, \m{p\geq 1} \ un entier, 
\m{n_0,\ldots,n_p} des entiers positifs, et pour \ \m{0\leq i\leq p}, 
\m{1\leq j\leq n_i}, \ \m{\ke^{(i)}_j} un faisceau coh\'erent sur $X$ et 
\m{M^{(i)}_j} un espace vectoriel non nul de dimension finie. On pose, pour
\m{0\leq i\leq p}
$$\ke_i \ = \ \som_{1\leq j\leq n_i}(\ke^{(i)}_j\ot M^{(i)}_j).$$
On suppose que les faisceaux \m{\ke^{(i)}_j} sont simples, et que
\ \m{\Hom(\ke^{(i)}_j, \ke^{(i')}_{j'}) \ = \ \nsp} \
si \ \m{i>i'}, ou \ \m{i=i'}, \m{j>j'}. Soit \m{\kw} la vari\'et\'e des
complexes
$$\ke_0\lra\ke_1\lra\ldots\lra\ke_p,$$
sur laquelle op\`ere le groupe alg\'ebrique
$$G \ = \ \Aut(\ke_0)\times\ldots\times\Aut(\ke_p).$$
Le sous-groupe unipotent $H$ est constitu\'e des \m{(g_0,\ldots,g_p)} tels que
pour tous $i$, $j$, la composante 
$$\ke^{(i)}_j\ot M^{(i)}_j\lra\ke^{(i)}_j\ot M^{(i)}_j$$
de \m{g_i} soit l'identit\'e. Le sous-groupe r\'eductif \m{G_{red}} est
constitu\'e des \m{(g_0,\ldots,g_p)} tels que pour tous $i$,$j$ ont ait
$$g_i(\ke^{(i)}_j\ot M^{(i)}_j)\ \subset \ \ke^{(i)}_j\ot M^{(i)}_j.$$
On a un isomorphisme
$$G_{red} \ \simeq \ \prod_{0\leq i\leq p,1\leq j\leq n_i} GL(M^{(i)}_j).$$
Si on veut retrouver une action sur une vari\'et\'e projective, il convient de
consid\'erer plut\^ot la vari\'et\'e projective \m{\P(W)} sur laquelle
op\`ere le groupe \m{G/\C^{*}}.

L'action de \m{G_{red}} sur \m{\kw} est un cas particulier des actions
\'etudi\'ees par A. King dans \cite{king}. Une lin\'earisation de l'action de
\m{G/\C^{*}} sur \m{\P(W)} est d\'efinie par une suite
\ \m{\Lambda=(\lambda_{ij})_{0\leq i\leq p,1\leq j\leq n_i}} \ de nombres
rationnels non nuls telle que
$$\sigg_{0\leq i\leq p,1\leq j\leq n_i}\lambda_{ij}\dim(M^{(i)}_j) \ = \ 0.$$
On appelle $\Lambda$ une {\em polarisation} de l'action de \m{G} sur \m{\kw}.
Un complexe
$$\ke_0 \ \hfl{f_0}{} \ \ke_1 \ \hfl{f_1}{} \ \ldots \ \hfl{f_{p-1}}{} \ \ke_p$$
est \m{G_{red}}-semi-stable (resp. \m{G_{red}}-stable) relativement \`a 
\m{\Lambda} si et seulement si pour tous sous-espaces vectoriels
$${M'}^{(i)}_j \ \subset M^{(i)}_j,$$
avec \ \m{({M'}^{(i)}_j)\not = (\nsp)} ou \m{( M^{(i)}_j)}, tels que
$$f_i(\som_{1\leq j\leq n_i}(\ke^{(i)}_j\ot {M'}^{(i)}_j)) \ \subset
\som_{1\leq j\leq n_{i+1}}(\ke^{(i+1)}_j\ot {M'}^{(i+1)}_j))$$
pour \ \m{0\leq i<p}, on a
$$\sigg_{0\leq i\leq p,1\leq j\leq n_i}\lambda_{ij}\dim({M'}^{(i)}_j) 
\ \leq \ 0 \ \ \ \ {\rm (resp. } \ < \ {\rm )}.$$

Le probl\`eme de la construction de quotients par $G$ d'ouverts de $\kw$
a \'et\'e abord\'e dans \cite{dr_tr} et \cite{dr2}, dans le cas \ \m{p=1}
(les complexes sont alors en fait des morphismes).  
On va donner ici une m\'ethode de construction de quotients par $G$ d'ouverts
$G$-invariants de $\kw$, qui est une g\'en\'eralisation de la m\'ethode
employ\'ee dans \cite{dr_tr}. 
On donnera en particulier des exemples de quotients
par $G$ de l'ouvert des points $G$-semi-stables. On ne peut pas obtenir en 
g\'en\'eral des {\em bons quotients}, mais ce qu'on appelle des {\em
quasi-bons quotients} (cf. chapitre 2). 

\vskip 1.5cm

\begin{sub}{\bf Mutations constructives}\end{sub}

Pour construire des bons quotients par $G$ d'ouverts $G$-invariants de $\kw$,
on introduit une nouvelle vari\'et\'e de complexes \m{\kw'}, sur laquelle le
groupe \m{G'} qui agit est r\'eductif. On tente ensuite d'\'etablir une relation
entre les quotients par \m{G'} d'ouverts \m{G'}-invariants de \m{\kw'} et les
quotients par $G$ d'ouverts $G$-invariants de $\kw$.

La m\'ethode est bas\'ee sur le r\'esultat suivant (cf. \para 3) : soient
\m{\ke}, $\kf$, $\kg$, $\Gamma$ des faisceaux coh\'erents sur $X$ et $M$ un
espace vectoriel de dimension finie. On suppose que le morphisme d'\'evaluation
$$\eva{\Gamma}{\kg}$$
est surjectif. Soit \m{\ke'} son noyau. On suppose aussi que la composition
$$\Hom(\ke,\Gamma)\ot\Hom(\Gamma,\kg)\lra\Hom(\ke,\kg)$$
est surjective. Soit 
$$(1) \ \ \ \ \ \ \ke \ \hfl{A}{} \ (\Gamma\ot M)\oplus\kg \ \hfl{B}{} \ \kf$$
un complexe. Alors on peut associer \`a \m{(1)} un complexe
$$(2) \ \ \ \ \ \
\ke\oplus\ke' \ \hfl{\alpha}{} \ \Gamma\ot N \ \hfl{\beta}{} \ \kf,$$
avec
\ \m{N \ = \ \Hom(\Gamma,\kg)\oplus M} \
et tel que 
$$\ker(\alpha)\simeq\ker(A), \ \ \ 
\ker(\beta)/\imm(\alpha)\simeq\ker(B)/\imm(A), \ \ \
\coker(\beta)\simeq\coker(\alpha).$$
La r\'eciproque est aussi vraie, si on part d'un complexe \m{(2)} tel que 
$\alpha$ induise une injection
$$\Hom(\ke',\Gamma)\lra N,$$
et les deux transformations sont inverses l'une de l'autre (\`a l'action
pr\`es des groupes d'automorphismes des complexes). Le passage de \m{(1)} \`a
\m{(2)} consiste \`a effectuer une {\em mutation \`a gauche} de $\kg$. Cette
notion a \'et\'e introduite dans l'\'etude des {\em \fexcs} (cf. \cite{dr_lp}, 
\cite{dr1}, \cite{go_ru}). On notera \m{\kw_0} (resp. \m{\kw'_0}) l'espace des
complexes \m{(1)} (resp. \m{(2)}), et \m{G_0} (resp. \m{G'_0}) le groupe
alg\'ebrique agissant sur \m{\kw_0} (resp. \m{\kw'_0}).

La transformation qui fait passer de \m{(1)} \`a \m{(2)} est purement
formelle (cf. \para 3.2). Le complexe \m{(2)} associ\'e \`a \m{(1)} n'est pas en
g\'en\'eral unique, mais sa \m{G'_0}-orbite l'est, et ne d\'epend que de la
\m{G_0}-orbite de \m{(1)}. On obtient ainsi (sous certaines hypoth\`eses) une
bijection
$$\kw_0/G_0 \ \simeq U_0/G'_0,$$
\m{U_0} d\'esignant l'ouvert de \m{\kw'_0} constitu\'e des complexes \m{(2)} tels que \m{\alpha} induise une injection \ \m{\Hom(\ke',\Gamma)\lra N}. Cette
correspondance est {\em alg\'ebrique} dans le sens suivant : si un ouvert
invariant d'un des espaces de complexes admet un quasi-bon quotient, l'ouvert
correspondant de l'autre c\^ot\'e admet aussi un quasi-bon quotient, et les deux
quotients sont isomorphes. La bijection pr\'ec\'edente est en fait un
{\em quasi-isomorphisme fort} (cf. chapitre 2). Cela entra\^ine que certaines
propri\'et\'es v\'erifi\'ees par un quasi-bon quotient d'un ouvert invariant
d'un des espaces de complexes seront automatiquement v\'erifi\'ees par le
quotient de l'ouvert correspondant de l'autre espace de complexes (cf. \para
2.3, concernant la descente sur les quotients de fibr\'es vectoriels). 
On a aussi une notion similaire de {\em mutation \`a droite}.

Pour appliquer ce qui pr\'ec\`ede aux vari\'et\'es de complexes de type
$$\ke_0\lra\ke_1\lra\ldots\lra\ke_p$$
on proc\`ede de la fa\c con suivante : on \'ecrit le dernier terme
$$\ke_p \ = (\Gamma\ot M)\oplus\kg,$$
avec
$$\Gamma=\ke^{(p)}_1, \ \ M=M^{(p)}_1, \ \
\kg=\som_{2\leq j\leq n_p}(\ke^{(p)}_j\ot M^{(p)}_j).$$
Dans ce cas, on a
$$\ke' \ = \ \som_{2\leq j\leq n_p}({\ke'}^{(p)}_j\ot M^{(p)}_j),$$
\m{{\ke'}^{(p)}_j} d\'esignant le noyau du morphisme d'\'evaluation
$$\eva{\ke^{(p)}_1}{\ke^{(p)}_j}$$
(suppos\'e surjectif). Les complexes obtenus par mutation sont du type
$$\ke_0\lra\ke_1\lra\ldots\lra\ke'_{p-1}\lra\ke^{(p)}_1\ot N,$$
avec
$$N \ = \ M^{(p)}_1\oplus\biggl(\som_{2\leq j\leq n_p}
(\Hom(\ke^{(p)}_1, \ke^{(p)}_j)\ot M^{(p)}_j)\biggr)$$
et
$$\ke'_{p-1} \ = \ \ke_{p-1}\oplus\biggl(\som_{2\leq j\leq n_p}
({\ke'}^{(p)}_j\ot M^{(p)}_j)\biggr).$$
On peut continuer en \'ecrivant
$$\ke'_{p-1} \ = (\Gamma\ot M)\oplus\kg,$$
avec
$$\Gamma=\ke^{(p-1)}_1, \ \ M=M^{(p-1)}_1, \ \
\kg=\biggl(\som_{2\leq j\leq n_{p-1}}(\ke^{(p-1)}_j\ot M^{(p-1)}_j)\biggr)\oplus
\biggl(\som_{2\leq j\leq n_p}({\ke'}^{(p)}_j\ot M^{(p)}_j)\biggr).$$
On peut ainsi proc\'eder \`a
$$q \ = \sigg_{0\leq i\leq p}n_i$$
mutations successives, et on obtient finalement un complexe du type
$$\kf_0\ot N_0\lra\ldots\lra\kf_q\ot N_q,$$
o\`u le groupe \m{G'} qui op\`ere est r\'eductif. Il existe un autre chemin
possible, en effectuant des mutations \`a droite en partant du terme de gauche. 

Supposons fix\'ee une polarisation \m{\Lambda} de l'action de $G$ sur $\kw$. On
peut alors d\'efinir naturellement une polarisation \m{\Lambda'} de l'action
de \m{G'} sur l'espace \m{\kw'} des complexes pr\'ec\'edents. Il reste 
\`a \'etudier les relations qu'il y a entre la \m{G}-(semi-)stabilit\'e des
complexes de \m{\kw} relativement \`a \m{\Lambda}, et la 
\m{G'}-(semi-)stabilit\'e
des complexes de \m{\kw'} relativement \`a \m{\Lambda'}. Il est toujours vrai
que si la mutation est \m{G'}-(semi-)stable relativement \`a \m{\Lambda'}, le
complexe d'origine est \m{G}-(semi-)stable relativement \`a \m{\Lambda}. La
r\'eciproque est vraie si on impose des conditions \`a \m{\Lambda}. Il faut 
ensuite montrer que tous les complexes \m{G'}-semi-stables de \m{\kw'} sont 
(\`a l'action de \m{G'} pr\`es) des mutations de complexes de $\kw$. Cela n'est
vrai que si on impose encore d'autres conditions \`a \m{\Lambda}. On obtient
alors l'existence d'un quasi-bon quotient projectif de l'ouvert des points 
$G$-semi-stables de $\kw$. En consid\'erant les mutations \`a droite, on
obtient g\'en\'eralement d'autres valeurs de $\Lambda$ pour lesquelles il
existe un quasi-bon quotient projectif.

\vskip 1.5cm

\begin{sub}{\bf Vari\'et\'es de modules de complexes}\end{sub}

On consid\`ere dans le \para 4 des complexes du type
$$(3) \ \ \ \ \ \
\ke_1\ot L_1\lra(\kf_1\ot M_1)\oplus(\kf_2\ot M_2)\lra\kg_1\ot N_1,$$
o\`u \m{\ke_1}, \m{\kf_1}, \m{\kf_2}, \m{\kg_1} sont des faisceaux coh\'erents
simples sur $X$, et \m{L_1}, \m{M_1}, \m{M_2}, \m{N_1} des espaces vectoriels de 
dimension finie. Quelques hypoth\`eses doivent \^etre faites, notamment que
le morphisme d'\'evaluation
$$\eva{\kf_1}{\kf_2}$$
est surjectif. On note \m{\kh_1} son noyau. En effectuant une premi\`ere
mutation \`a gauche on associe au complexe \m{(3)} un complexe
$$(4) \ \ \ \ \ \
(\ke_1\ot L_1)\oplus(\kh_1\ot M_2)\lra\kf_1\ot P_1\lra\kg_1\ot N_1,$$
avec
$$P_1 \ = \ (\Hom(\kf_1,\kf_2)\ot M_2)\oplus M_1.$$
On suppose ensuite que le morphisme d'\'evaluation
$$\eva{\ke_1}{\kh_1}$$
est surjectif. Soit \m{\kk_1} son noyau. En effectuant une seconde mutation 
on obtient un complexe
$$(5) \ \ \ \ \ \
\kk_1\ot M_2\lra\ke_1\ot Q_1\lra\kf_1\ot P_1\lra\kg_1\ot N_1,$$
avec
$$Q_1 \ = \ (\Hom(\ke_1,\kh_1)\ot M_2)\oplus L_1,$$
et ici le groupe qui op\`ere est r\'eductif. 

On en d\'eduit dans le th\'eor\`eme 4.4 l'existence de quasi-bons quotients 
projectifs d'ouverts de complexes $G$-semi-stables de type \m{(3)} (pour
certaines polarisations). 

Dans le \para 4.5 on \'etudie le cas des complexes
$$\ko(-2)\lra(\ko(-1)\ot M_1)\oplus\ko\lra\ko(1)$$
sur \m{\P_n}, \m{M_1} \'etant un espace vectoriel tel que \
\m{0<\dim(M_1)<n+1}. Dans le cas de \m{\P_2} et \ \m{\dim(M_1)=3}, on
obtient trois types de quotients distincts, dont deux lisses. Un des quotients
lisses peut \^etre obtenu de mani\`ere \'el\'ementaire. Le second est non
trivial.

\vskip 1.5cm

\begin{sub}{\bf Vari\'et\'es de modules de morphismes}\end{sub}

Dans le \para 5 on consid\`ere des morphismes du type
$$(\ke_1\ot M_1)\oplus(\ke_2\ot M_2)\lra\kf_1\ot N_1,$$
o\`u \m{\ke_1}, \m{\ke_2}, \m{\kf_1} sont des faisceaux coh\'erents simples sur
$X$, et \m{M_1}, \m{M_2}, \m{N_1} des espaces vectoriels de dimension finie.
On rappelle dans les \para 5.1 et 5.2 la construction de vari\'et\'es de
modules de morphismes $G$-semi-stables de ce type (pour certaines polarisations)
effectu\'ee dans \cite{dr_tr}. On proc\`ede simplement ici \`a une seule
mutation \`a gauche pour obtenir une action d'un groupe r\'eductif. Les
quotients obtenus par cette m\'ethode sont des bons quotients projectifs, 
et les ouverts 
correspondant aux morphismes $G$-stables sont des quotients g\'eom\'etriques.

Dans le \para 5.3 on emploie des mutations \`a droite. Il faut alors deux
mutations successives pour obtenir une action d'un groupe r\'eductif. On
obtient des bons quotients de l'ouvert des morphismes $G$-semi-stables,
pour
d'autres polarisations qu'avec la m\'ethode pr\'ec\'edente. 

Dans la \para 5.4 on donne des exemples de constructions de vari\'et\'es de
modules de morphismes au moyen des m\'ethodes de \cite{dr_tr}, \cite{dr2}
et du \para 5.3. On donne un cas o\`u il n'y a pas de quotient
g\'eom\'etrique de l'ouvert des points stables.

\vskip 1.5cm

\begin{sub}{\bf Mutations non constructives}\end{sub}

D'autres sortes de mutations ont \'et\'e d\'efinies dans \cite{dr2}. On 
pourrait les appeler des {\em mutations non constructives}. Dans \cite{dr2}
elles sont appliqu\'ees \`a des morphismes, mais on peut sans difficult\'e 
\'etendre leur d\'efinition aux complexes. Elles peuvent aussi servir \`a
construire des bons quotients d'ouverts de points $G$-semi-stables, pour
d'autres polarisations que celles qui sont accessibles par les m\'ethodes
d\'ecrites ici. La d\'efinition des mutations non constructives est bas\'ee sur
le r\'esultat suivant : soient \m{\ke}, \m{\ke'}, \m{\Gamma}, \m{\kg} et \m{\kf} 
des faisceaux coh\'erents sur $X$. On suppose que le morphisme canonique
$$\coeva{\ke'}{\Gamma}$$
est injectif et on note \m{\kf_0} son conoyau. On suppose aussi que
$$\Ext^1(\kf_0,\kf)=\Ext^1(\ke,\Gamma)=\Ext^1(\Gamma,\kf)=\nsp.$$
Soient $M$ un espace vectoriel de dimension finie et
$$\ke\oplus\ke' \ \hfl{A}{} \ (\Gamma\ot M)\oplus\kf \ \hfl{B}{} \ \kg$$
un complexe tel que l'application lin\'eaire 
$$\lambda : \Hom(\ke',\Gamma)^*\lra M$$
d\'eduite de $A$ soit surjective. Alors il existe un complexe
$$\ke\oplus(\Gamma\ot\ker(\lambda)) \ \hfl{\alpha}{} \ \kf_0\oplus\kf \
\hfl{\beta}{} \ \kg$$
tel que 
$$\ker(\alpha)\simeq\ker(A), \ \ \ 
\ker(\beta)/\imm(\alpha)\simeq\ker(B)/\imm(A), \ \ \
\coker(\beta)\simeq\coker(B).$$
Ce r\'esultat admet aussi une r\'eciproque. 

Le diff\'erence essentielle entre les mutations constructives et les 
mutations non constructives est la suivante : dans le premier cas, on
effectue la mutation d'une paire de faisceaux situ\'es dans le m\^eme
terme du complexe, et dans le second cas on effectue la mutation d'une paire 
de faisceaux situ\'es dans des termes adjacents du complexe.

\vskip 2.5cm

\section{Quotients alg\'ebriques}

\begin{sub}{\bf Quasi-bons quotients}\end{sub}

Soit $G$ un groupe alg\'ebrique. On appelle {\em $G$-espace} une vari\'et\'e
alg\'ebrique $X$ munie d'une action alg\'ebrique de $G$. Rappelons qu'on
appelle {\em bon quotient} de $X$ par $G$ un morphisme
$$\pi : X\lra M$$
(o\`u $M$ est une vari\'et\'e alg\'ebrique) tel que : 

(i) Le morphisme $\pi$ est $G$-invariant, affine et surjectif.

(ii) Si $U$ est un ouvert de $M$, alors on a
\ \m{\ko(U) \ \simeq \ \ko(\pi^{-1}(U))^G}.

(iii) Si $F_1$, $F_2$ sont des sous-vari\'et\'es ferm\'ees
$G$-invariantes disjointes de $X$, alors \m{\pi(F_1)} et \m{\pi(F_2)} sont des
sous-vari\'et\'es disjointes de $M$. 

\medskip

Cette d\'efinition est particuli\`erement
bien adapt\'ee \`a l'\'etude des actions de groupes alg\'ebriques r\'eductifs,
car on sait que dans ce cas il existe toujours un bon quotient si $X$ est
affine. On utilisera une notion l\'eg\`erement diff\'erente :

\bigskip

\begin{defin}
On appelle {\em quasi-bon quotient} de $X$ par $G$ un morphisme
$$\pi : X\lra M$$
(o\`u $M$ est une vari\'et\'e alg\'ebrique) qui est $G$-invariant,
surjectif et tel que les conditions (ii) et (iii) pr\'ec\'edentes soient 
v\'erifi\'ees.
\end{defin}

\bigskip

On dit parfois par abus de language que $M$ est le quasi-bon quotient de $X$
par $G$. On note comme dans le cas des bons quotients
\ \m{M \ = \ X//G}.
Il est clair qu'un bon quotient est un quasi-bon quotient.

\vskip 1.5cm

\begin{sub}{\bf Quasi-isomorphismes}\end{sub}

Soient $G$, $G'$ des groupes alg\'ebriques, $X$ un $G$-espace et $X'$ un
\m{G'}-espace. 

\bigskip

\begin{defin}
1 - On appelle {\em quasi-morphisme} de $X$ vers $X'$ une application
$$\phi : X/G\lra X'/G'$$
telle que pour tout point $x$ de $X$ il existe un ouvert de Zariski $U$ de $X$
contenant $x$ et un morphisme \ \m{U\lra X'} \ induisant $\phi$.

\medskip

2 - On appelle {\em quasi-morphisme fort} de $X$ vers $X'$ la
donn\'ee d'un {\em quasi-morphisme}
$$\phi : X/G\lra X'/G',$$
d'un recouvrement ouvert \m{(U_i)_{i\in I}} de X, et d'une famille
\m{(\phi_i)_{i\in I}} de rel\`evements de \m{\phi} ,
\ \m{\phi_i : U_i\lra X'},
telle que :

(i) pour tout \m{j\in I} et \m{g\in G}, le morphisme
$$gU_j\lra X'$$
$$x\longmapsto \phi_j(g^{-1}x)$$
appartient \`a la famille \m{(\phi_i)}.

(ii) Pour tous \ \m{i,j\in I}, et tout \ \m{x\in U_i\cap U_j} \
il existe un voisinage $V$ de $x$ dans \m{U_i\cap U_j} et un morphisme \ 
\m{\lambda_{ij} : V\lra G'} \ tel que \
\m{{\phi_j}_{\mid V} \ = \ \lambda_{ij}{\phi_i}_{\mid V}},

(iii) Pour tout \ \m{i\in I} \ et \ \m{x\in U_i}, il existe un
voisinage $V$ de \m{(e,x)} dans \ \m{G\times U_i} \ et un morphisme 
\ \m{\gamma : V\lra G'} \
tel que \ \m{\gamma(e,x)=e} \  et que pour tous \ \m{(g,y)\in V}, on
ait \ \m{gy\in U_i} \ et \
\m{\phi_i(gy)=\gamma(g,y)\phi_i(y)}.
\end{defin}

\bigskip

On appelle {\em quasi-isomorphisme} de $X$  vers $X'$ une bijection
$$\phi : X/G\lra X'/G'$$
qui est un quasi-morphisme ainsi que son inverse. 

\bigskip

\begin{defin}
Soit \ \m{\sigma = (\phi : X/G\lra X'/G', (\phi_i : U_i\lra X'))} \ un
quasi-morphisme fort. On appelle {\em carte} de \m{\sigma} un rel\`evement
local
\ \m{f : U\lra X'} de \m{\phi}
($U$ \'etant un ouvert non vide de $X$) tel que pour tout 
\m{x\in U} les propri\'et\'es suivantes soient v\'erifi\'ees :

\noindent (i) Pour tout \m{i\in I} il existe un voisinage $V$ de $x$ dans 
$U\cap U_i$ et un morphisme
\m{\lambda_{i} : V\lra G'} \ tel que \
\m{{f}_{\mid V} \ = \ \lambda_{i}{\phi_i}_{\mid V}},

\noindent (ii) Il existe un voisinage $V$ de
\m{(e,x)} dans \ \m{G\times U} \ et un morphisme 
\ \m{\gamma : V\lra G'} \
tel que \ \m{\gamma(e,x)=e} \  et que pour tous \ \m{(g,y)\in V}, on
ait \ \m{gy\in U} \ et
\ \m{\phi_i(gy)=\gamma(g,y)\phi(y)}.
\end{defin}

\bigskip

On dit que deux quasi-morphismes forts de $X$ dans \m{X'} sont {\em
\'equivalents} si les cartes de l'un sont aussi des cartes de l'autre. 

\bigskip

Il est clair que la composition de deux quasi-morphismes en est un. On va
d\'efinir ce qu'on entend par {\em composition} de deux quasi-morphismes forts.
Soient $G$, \m{G'}, \m{G''} des groupes alg\'ebriques, \m{X} un $G$-espace,
\m{X'} un \m{G'}-espace et \m{X''} un \m{G''}-espace. Soient 
$$\sigma = (\phi : X/G\lra X'/G', (\phi_i : U_i\lra X')_{i\in I}),$$
$$\sigma' = (\phi' : X'/G\lra X''/G'', (\phi'_j : U'_j\lra X'')_{j\in J})$$
des quasi-morphismes forts. Soient
$$\psi=\phi'\circ\phi : X/G\lra X''/G'',$$
et pour \m{i\in I}, \m{j\in J}, \ \m{V_{ij}=U_i\cap\phi^{-1}(U'_j)}, et
$$\psi_{ij} = \phi'_j\circ\phi_i : V_{ij}\lra X''.$$
On pose
$$\tau = (\psi : X/G\lra X''/G'', 
(\psi_{i,j} : V_{ij}\lra X'')_{(i,j)\in I\times J}).$$

\bigskip

\begin{xlemm}
$\tau$ est un quasi-morphisme fort de $X$ dans $X'$.
\end{xlemm}

\bigskip

\begin{proof} Soit \ \m{x\in V_{ij}}. Pour \m{(g,y)} dans un voisinage 
convenable de \m{(e,x)}, on a
$$\psi_{ij}(gy) = \phi'_j\circ\phi_i(gy)=\phi'_j(\gamma_i(g,y)\phi_i(y))
=\gamma'_j(\gamma_i(g,y),\phi_i(y))\psi_{ij}(y),$$
\m{\gamma_i} (resp. \m{\gamma'_j}) \'etant un morphisme \`a valeurs dans
\m{G'} (resp. \m{G''}) d\'efini sur un voisinage de \m{(e,x)} (resp.
\m{(e,\phi_i(x))}). En posant
$$\gamma_{ij}(g,y) = \gamma'_j(\gamma_i(g,y),\phi_i(y))$$
on obtient donc \
\m{\psi_{ij}(gy) = \gamma_{ij}(g,y)\psi_{ij}(y)}.
D'autre part, soient \ \m{i,i'\in I},\m{j,j'\in J} \ et \hfil\break
\m{x\in V_{i,j}\cap V_{i',j'}}. Avec des notations \'evidentes, on a, pour $y$
dans un voisinage de $x$
$$\phi'_{j'}\circ\phi_{i'}(y)=\phi'_{j'}(\lambda_{ii'}(y)\phi_i(y)).$$
Posons \
\m{g'_0=\lambda_{ii'}(x)\in G'}.
Alors on a
$$\phi'_{j'}(\lambda_{ii'}(y)\phi_i(y))=\gamma'_{j'}(\lambda_{ii'}(y){g'_0}^{-1},
g'_0\phi_i(y))\phi'_{j'}(g_0\phi_i(y))=\gamma'_{j'}(\lambda_{ii'}(y){g'_0}^{-1},
g'_0\phi_i(y))\phi'_{k}(\phi_i(y))$$
pour un $k$ convenable dans $J$. On obtient finalement
$$\phi'_{j'}\circ\phi_{i'}(y)=\gamma'_{j'}(\lambda_{ii'}(y){g'_0}^{-1},
g'_0\phi_i(y))\lambda'_{kj}(\phi_i(y))\phi'_j\circ\phi_i(y),$$
c'est-\`a-dire
$$\phi'_{j'}\circ\phi_{i'}(y)=\theta_{i'j',ij}(y)\phi'_{j}\circ\phi_{i}(y),$$
avec
$$\theta_{i'j',ij}(y)=\gamma'_{j'}(\lambda_{ii'}(y){g'_0}^{-1},
g'_0\phi_i(y))\lambda'_{kj}(\phi_i(y)).$$ \end{proof}

\bigskip

On d\'efinit de mani\`ere \'evidente le quasi-morphisme fort {\em identit\'e}
\m{I_X} de $X$ dans $X$. Il est clair que la composition \`a droite ou \`a 
gauche d'un quasi-morphisme fort avec \m{I_X} donne un quasi-morphisme fort
\'equivalent. On d\'efinit un {\em quasi-isomorphisme fort} $\sigma$ de $X$ 
dans \m{X'} comme \'etant un quasi-morphisme fort de $X$ dans \m{X'} tel qu'il 
existe un quasi-morphisme fort \m{\sigma'} de \m{X'} dans $X$ tel que
\m{\sigma\circ\sigma'} (resp. \m{\sigma'\circ\sigma}) soit \'equivalent \`a
\m{I_{X'}} (resp. \m{I_X}).

\bigskip

{\bf Exemples : } 1 - Supposons qu'il existe un morphisme 
\ \m{\phi : X\lra X'} \ 
compatible avec un morphisme de groupes alg\'ebriques 
\ \m{\alpha : G\lra G'},
c'est-\`a-dire que pour tous \m{g\in G}, \m{x\in X}, on a \ 
\m{\phi(gx)=\alpha(x)\phi(x)}. Soit 
$$\ov{\phi} : X/G\lra X'/G'$$
l'application d\'eduite de \m{\phi}. Alors \m{(\ov{\phi}, \phi)} est un
quasi-morphisme fort de $X$ dans \m{X'}. C'est l'exemple le plus simple. 

\medskip

2 - Supposons qu'il existe un quasi-isomorphisme fort \
\m{\sigma : X\lra X'}, et soit $U$ un ouvert $G$-invariant de $X$. 
Soit \m{U'} l'ouvert
\m{G'}-invariant correspondant de \m{X'}. Alors \m{\sigma} induit un
quasi-isomorphisme fort \ \m{U\lra U'}.

\bigskip

Il est clair que si 
$$\phi : X/G\lra X'/G'$$ 
est un quasi-isomorphisme (non n\'ecessairement fort), 
\m{\phi} induit une bijection entre l'ensemble des ouverts (resp. ferm\'es) 
$G$-invariants de $X$ et
l'ensemble des ouverts (resp. ferm\'es) \m{G'}-invariants de \m{X'}. Si $U$
est un ouvert $G$-invariant de $X$, et \m{U'} l'ouvert \m{G'}-invariant
correspondant de \m{X'}, \m{\phi} induit un isomorphisme d'anneaux
$$\ko(U)^G \ \simeq \ \ko(U')^{G'}.$$
Plus g\'en\'eralement, si $Y$ est une vari\'et\'e alg\'ebrique, les morphismes
$G$-invariants \ \m{X\lra Y} \ s'identifient de mani\`ere \'evidente aux
morphismes \m{G'}-invariants \ \m{X'\lra Y}.

\bigskip

\begin{xprop}
Si $X$ et $X'$ sont quasi-isomorphes, il existe un quasi-bon quotient de
$X$ par $G$ si et seulement si il existe un quasi-bon quotient de \m{X'} par
\m{G'}, et dans ce cas les deux quotients sont isomorphes.
\end{xprop}

\bigskip

\begin{proof} Supposons qu'il existe un quasi-bon quotient 
$$\pi : X\lra M$$
de $X$ par $G$. Soit
$$\phi : X/G\lra X'/G'$$
un quasi-isomorphisme. Puisque \m{\pi} est $G$-invariant, il d\'efinit un
morphisme \m{G'}-invariant
$$\pi' : X'\lra M.$$ 
Il est imm\'ediat que c'est un quasi-bon quotient de \m{X'} par \m{G'}.
\end{proof}

\bigskip

{\bf Remarque : } La proposition 2.2 est vraie pour des bons quotients
si on suppose en plus que les vari\'et\'es $X$ et \m{X'} sont affines et si
tout ouvert affine d'une de ces vari\'et\'es est le compl\'ementaire d'une
hypersurface (c'est le cas par exemple lorsque $X$ et \m{X'} sont
{\em factorielles}, c'est-\`a-dire que leurs anneaux de fonctions r\'eguli\`eres
sont factoriels). En effet dans ce cas les ouverts affines $G$-invariants 
(resp. \m{G'}-invariants) de $X$ (resp. \m{X'}) sont alors les 
compl\'ementaires des hypersurfaces $G$-invariantes
(resp. \m{G'}-invariantes), et il est ais\'e de voir que via \m{\phi} les
hypersurfaces invariantes de $X$ correspondent exactement \`a celles de
\m{X'}. C'est ce qui se produit dans \cite{dr2}, o\`u $X$ et \m{X'} sont des
ouverts d'espaces affines.

\vskip 1.5cm

\begin{sub}{\bf Lemme de descente}\end{sub}

On montre ici que la notion de quasi-isomorphisme fort pr\'eserve une 
propri\'et\'e importante : la 
possibilit\'e de descendre au quotient des $G$-fibr\'es ad\'equats sur $X$.
Rappelons qu'un {\em $G$-fibr\'e vectoriel} sur $X$ est un fibr\'e vectoriel
alg\'ebrique $F$ sur $X$ muni d'une action alg\'ebrique lin\'eaire de $G$, au
dessus de l'action de $G$ sur $X$. S'il existe un quasi-bon quotient
\ \m{\pi : X\lra M}, on dit que $F$ {\em descend} \`a $M$ s'il existe un fibr\'e
vectoriel alg\'ebrique $E$ sur $M$ et un isomorphisme de $G$-fibr\'es
$$F \ \simeq \ \pi^*(E).$$

\bigskip

\begin{defin}
On dit qu'un $G$-fibr\'e vectoriel $F$ sur $X$ est {\em admissible} si pour
tout point $x$ de $X$, le stabilisateur de $x$ dans $G$ agit trivialement sur
$F_x$.
\end{defin}

\bigskip

Il est clair que s'il existe un quasi-bon quotient \ \m{M=X//G}, et si $F$
descend \`a $M$, alors $F$ est admissible. 
On d\'emontre dans \cite{dr_na} le r\'esultat suivant :

\bigskip

\begin{xlemm}
{\em (Lemme de descente)} Si le groupe $G$ est r\'eductif et s'il existe un bon
quotient \ \m{\pi : X\lra M}, tout $G$-fibr\'e admissible sur $X$ descend \`a
$M$.
\end{xlemm}

\bigskip

Une justification de la notion de quasi-isomorphisme fort est le r\'esultat
suivant :

\bigskip

\begin{xprop}
On suppose donn\'e un quasi-isomorphisme fort de $X$ dans \m{X'}. Soit $r$
un entier. Alors

1 - Il existe une bijection canonique entre l'ensemble des classes
d'isomorphisme de $G$-fibr\'es vectoriels admissibles de rang $r$ sur $X$ et
l'ensemble des classes d'isomorphisme de \m{G'}-fibr\'es vectoriels 
admissibles de rang $r$ sur \m{X'}.

2 - On suppose qu'il existe un quasi-bon quotient 
$$M \ = \ X//G \ = \ X'//G'.$$
Alors un $G$-fibr\'e vectoriel admissible sur $X$ descend \`a $M$ si et
seulement si le \m{G'}--fibr\'e vectoriel admissible correspondant sur \m{X'}
descend \`a $M$, et les fibr\'es vectoriels associ\'es sur $M$ sont les 
m\^emes. 
\end{xprop}

\bigskip

\begin{proof}
Soit \ \m{(\phi : X/G\lra X'/G', (\phi_i : U_i\lra X'))} \ le quasi-isomorphisme
fort.
Soit $F'$ un $G'$-fibr\'e vectoriel admissible sur \m{X'}. On va lui associer
un $G$-fibr\'e vectoriel admissible $F$ sur $X$. On commence par construire
$F$ sans sa structure alg\'ebrique. Soit \ \m{\epsilon : X\lra X'} \ une
application au dessus de \m{\phi}. On pose
\ \m{F_\epsilon = \epsilon^*(F')}.
Si \ \m{\epsilon' : X\lra X'} \ est une autre application au dessus de \m{\phi},
il existe une application
\ \m{\theta : X\lra G'} \
telle que \ \m{\epsilon'=\theta.\epsilon}. Puisque $F'$ est un $G'$-fibr\'e
on obtient un isomorphisme
\ \m{F_\epsilon \simeq F_{\epsilon'}} \
qui en \m{x\in X} est la multiplication par \m{\theta(x)} 
\ \m{F'_{\epsilon(x)}\lra F'_{\epsilon'(x)}}.
Puisque $F'$ est admissible, l'isomorphisme pr\'ec\'edent est ind\'ependant du
choix de \m{\epsilon}. On peut donc d\'efinir sans ambiguit\'e 
\ \m{F \ = \ F_\epsilon}.

D\'efinissons maintenant l'action de $G$ sur $F$. Soient \m{g\in G}, \m{x\in X}.
Alors \m{\epsilon(x)} et \m{\epsilon(gx)} sont dans la m\^eme \m{G'}-orbite,
donc on a un isomorphisme canonique
$$F_x=F'_{\epsilon(x)}\lra F'_{\epsilon(gx)}=F_{gx},$$
qui ne d\'epend pas de $\epsilon$. Il est clair qu'on d\'efinit ainsi une
action de $G$ sur $F$, et que le stabilisateur de $x$ agit trivialement sur
\m{F_x}. 

On d\'efinit maintenant la structure alg\'ebrique sur $F$. Soient \m{i\in I}
et \m{x\in U_i}. Soit $V$ un voisinage de \m{\phi_i(x)} tel qu'on ait une
trivialisation
$$F'_{\mid V} \ \simeq \ko_V\ot\C^{r}.$$
On d\'efinit sur \ \m{U=U_i\cap\phi_i^{-1}(V)} \ une trivialisation locale 
$$F_U \ \simeq \phi^*(F'_{\mid V}) \simeq \ko_U\ot\C^{r}.$$
Il faut v\'erifier que ces trivialisations se recollent alg\'ebriquement. Cela
d\'ecoule imm\'ediatement de la condition (ii) de la d\'efinition d'un
quasi-morphisme fort. Le fait que l'action de $G$ est alg\'ebrique d\'ecoule
de la condition (iii).

Les autres assertions sont imm\'ediates et laiss\'ees au lecteur.
\end{proof}

\bigskip

Si $G$ est r\'eductif et s'il existe un bon quotient \ \m{M=X//G}, 
on en d\'eduit
que tout \m{G'}-espace fortement quasi-isomorphe \`a $X$ v\'erifie le
lemme de descente, c'est-\`a-dire que tout $G'$-fibr\'e vectoriel admissible
sur $X'$ descend \`a $M$.

\vskip 2.5cm

\section{Mutations constructives}

\begin{sub}{\bf Mutations constructives en termes de faisceaux}\end{sub}

\begin{subsub}{Un exemple simple}\end{subsub}

Soient \m{\ke}, \m{\kf}, \m{\kg}, \m{\Gamma} des faisceaux coh\'erents sur une
vari\'et\'e projective \m{X}, \m{M} un espace vectoriel de dimension 
finie non nul. On suppose que le morphisme d'\'evaluation
$$\Gamma\ot\Hom(\Gamma,\kg)\lra\kg$$
est surjectif. Soit \m{\ke'} son noyau. On suppose aussi que l'application
lin\'eaire canonique
$$\Hom(\Gamma,\kg)^*\lra\Hom(\ke',\Gamma)$$
est bijective, et que la composition
$$c : \Hom(\ke,\Gamma)\ot\Hom(\Gamma,\kg)\lra\Hom(\ke,\kg)$$
est surjective.

\bigskip

\begin{xprop} 
Soit
$$\ke\ \hfl{A}{}\ (\Gamma\ot M)\oplus\kg\ \hfl{B}{}\ \kf$$
un complexe, avec \m{A} injectif et \m{B} surjectif. Alors il existe un
complexe
$$\ke\oplus\ke'\ \hfl{\alpha}{}\ \Gamma\ot N\ \hfl{\beta}{}\ \kf,$$
avec \
\m{N=\Hom(\Gamma,\kg)\oplus M},
\m{\alpha} \'etant injectif, \m{\beta} surjectif, et
\ \m{\ker(\beta)/\imm(\alpha) \ \simeq \ \ker(B)/\imm(A)}. 
\end{xprop}

\medskip

\begin{proof} Le morphisme
$$\alpha : \ke\oplus\ke'\lra\Gamma\ot (M\oplus\Hom(\Gamma,\kg))$$
est la somme d'un morphisme 
$$\ke\lra\Gamma\ot (M\oplus\Hom(\Gamma,\kg))$$
provenant de \m{A} (qui existe car \m{c} est surjective), et de l'inclusion
\ \m{\ke'\lra\Gamma\ot\Hom(\Gamma,\kg)}.
On a un diagramme commutatif avec colonnes exactes :

\bigskip

\begin{picture}(360,230)
\put(135,220){$0$}
\put(280,220){$0$}
\put(133,170){$\ke'$}
\put(280,170){$\ke'$}
\put(120,120){$\ke'\oplus\ke$}
\put(230,120){$\Gamma\ot(M\oplus\Hom(\Gamma,\kg))$}
\put(135,70){$\ke$}
\put(260,70){$(\Gamma\ot M)\oplus\kg$}
\put(135,20){$0$}
\put(280,20){$0$}
\put(187,126){$\alpha$}
\put(195,76){$A$}

\put(137,214){\vector(0,-1){30}}
\put(137,164){\vector(0,-1){30}}
\put(137,114){\vector(0,-1){30}}
\put(137,64){\vector(0,-1){25}}
\put(282,214){\vector(0,-1){30}}
\put(282,164){\vector(0,-1){30}}
\put(282,114){\vector(0,-1){30}}
\put(282,64){\vector(0,-1){30}}

\put(145,172){\line(1,0){129}}
\put(145,174){\line(1,0){129}}
\put(157,123){\vector(1,0){70}}
\put(145,73){\vector(1,0){112}}
\end{picture}

\bigskip

On en d\'eduit que \m{\alpha} est injectif et que \ \m{\coker(\alpha) \ 
\simeq \coker(A)}. On d\'efinit le morphisme \m{\beta} par le carr\'e 
commutatif

\bigskip

\begin{picture}(360,90)
\put(105,80){$\Gamma\ot(M\oplus\Hom(\Gamma,\kg))$}
\put(280,80){$\kf$}
\put(120,30){$(\Gamma\ot M)\oplus\kg$}
\put(280,30){$\kf$}
\put(242,86){$\beta$}
\put(230,36){$B$}

\put(223,83){\vector(1,0){47}}
\put(200,33){\vector(1,0){70}}
\put(145,74){\vector(0,-1){30}}
\put(284,74){\line(0,-1){30}}
\put(286,74){\line(0,-1){30}}
\end{picture}

\bigskip

On en d\'eduit que \m{\beta} est surjectif, \m{\beta\circ\alpha=0} \ et
\ \ \m{\ker(\beta)/\imm(\alpha) \ \simeq \ \ker(B)/\imm(A)}.
\end{proof}

\bigskip

On a bien s\^ur une r\'eciproque :

\bigskip

\begin{xprop}
Soit
$$\ke\oplus\ke'\ \hfl{\alpha}{}\ \Gamma\ot N\ \hfl{\beta}{}\ \kf$$
un complexe, avec \m{\alpha} injectif et \m{\beta} surjectif. On suppose
que l'application lin\'eaire 
$$\lambda : \Hom(\ke',\Gamma)\lra N$$
d\'eduite de \m{\alpha} est injective. Alors il existe un complexe
$$\ke\ \hfl{A}{}\ (\Gamma\ot M)\oplus\kg\ \hfl{B}{}\ \kf$$
avec \ \m{M=\coker(\lambda)}, \m{A} \'etant injectif, \m{B} surjectif et
\ \m{\ker(B)/\imm(A) \ \simeq \ \ker(\beta)/\imm(\alpha)}.
\end{xprop}

\begin{proof} Analogue \`a la proposition pr\'ec\'edente. \end{proof}

\vskip 1cm

\begin{subsub}{Le cas g\'en\'eral}\end{subsub}

On se place dans la situation du \para 3.1. Soient \m{\ku}, \m{\kv}, \m{\kg_0}
des faisceaux coh\'erents sur \m{X}. On d\'emontre comme la proposition 3.1
la

\bigskip

\begin{xprop}
1 - Soit
$$0\lra\ku\lra\ke\ \hfl{A}{}\ (\Gamma\ot M)\oplus\kg\oplus\kg_0\ \hfl{B}{}\ 
\kf\lra\kv\lra 0$$
un complexe, exact en \m{\ku}, \m{\ke}, \m{\kf}, \m{\kv}. Alors il existe un
complexe
$$0\lra\ku\lra\ke\oplus\ke'\ \hfl{\alpha}{}\ (\Gamma\ot N)\oplus\kg_0\ 
\hfl{\beta}{}\ \kf\lra\kv\lra 0,$$
avec 
$$N = \Hom(\Gamma,\kg)\oplus M,$$
exact sauf au plus en \ \m{(\Gamma\ot N)\oplus\kg_0}, et tel que \
\m{\ker(\beta)/\imm(\alpha) \ \simeq \ \ker(B)/\imm(A)}. 

\medskip

2 - R\'eciproquement, si \m{N} est un espace vectoriel, et si on
a un complexe du second type exact en \m{\ku}, \m{\ke\oplus\ke'}, \m{\kf},
\m{\kv}, tel que que l'application lin\'eaire 
$$\lambda : \Hom(\ke',\Gamma)\lra N$$
d\'eduite de \m{\alpha} soit injective. Alors il existe un complexe du 
premier type, avec \
\m{M=\coker(\lambda)}, 
exact en \ \m{\ku}, \m{\ke}, \m{\kf}, \m{\kv}, et tel que
\ \ \m{\ker(B)/\imm(A) \ \simeq \ \ker(\beta)/\imm(\alpha)}.
\end{xprop}

\vskip 1.5cm

\begin{sub}{\bf Mutations constructives abstraites}\end{sub}

On d\'ecrit ici de mani\`ere abstraite la situation de la proposition 3.1
(sans tenir compte de l'injectivit\'e de \m{A} et \m{\alpha} et de la
surjectivit\'e de \m{B} et \m{\beta}). Il est possible de faire la
m\^eme chose dans le cas plus g\'en\'eral de la proposition 3.3. 
On \'etudie l'action de certains
groupes d'automorphismes sur l'espace de tous les complexes, et on
\'etudie la relation entre les orbites des deux types de complexes. 
Des hypoth\`eses suppl\'ementaires sont faites dans cette version 
abstraite, par exemple on suppose que \m{\Gamma} est simple. 

\vskip 1cm

\begin{subsub}{Espaces de complexes de type 1 (version simplifi\'ee)}
\end{subsub}

\noindent{\bf \ 3.2.1.1 }{\it D\'efinition}

\medskip

Soient \m{Z_1}, \m{Z_2}, \m{Z_3}, \m{Z_4} , \m{H}, \m{T}, \m{M} des espaces 
vectoriels de dimension finie. On pose
$$W_C = (Z_1\ot M)\oplus Z_2\oplus (Z_3\ot M^*)\oplus Z_4.$$
Soient
$$\sigma : Z_1\ot H\lra Z_2,$$
$$\sigma' : H\ot Z_4\lra Z_3,$$
$$\tau : Z_1\ot Z_3\lra T,$$
$$\tau' : Z_2\ot Z_4\lra T$$
des applications lin\'eaires. On suppose que \m{\sigma} est surjective, et
que \m{\sigma'} induit une inclusion \ \m{Z_4\subset H^*\ot Z_3}. On 
suppose aussi que le diagramme suivant est commutatif :
$$(D) \ \ \ \ \ \ \begin{CD}
Z_1\ot H\ot Z_4  @>\sigma\ot I_{Z_4}>>  Z_2\ot Z_4\cr
@V{I_{Z_1}\ot\sigma'}VV  @VV{\tau'}V\cr
Z_1\ot Z_3 @>>{\tau}> T\cr\end{CD}$$
Soit \ \m{Q_C\subset W_C} \ l'ensemble des points 
\m{(\phi_1,z_2,\phi_3,z_4)} tels que
$$\tau(\pline{\phi_1,\phi_3})+\tau'(z_2\ot z_4)=0,$$
\m{\pline{\phi_1,\phi_3}} d\'esignant l'image de \m{\phi_1\ot\phi_3} par la
contraction de \m{M}
$$Z_1\ot M\ot Z_3\ot M^*\lra Z_1\ot Z_3.$$

Soient \m{G_L}, $G_R$ et $G_0$ des groupes. On suppose que :
\begin{itemize}
\item[] $G_L$ op\`ere lin\'eairement \`a droite sur $Z_1$, $Z_2$,
$T$.
\item[] $G_R$ op\`ere lin\'eairement \`a gauche sur $Z_3$, $Z_4$,
$T$.
\item[] $G_0$ op\`ere lin\'eairement \`a gauche sur $Z_2$, $H$, et
lin\'eairement \`a droite sur $Z_4$.
\end{itemize}
On suppose que ces actions sont compatibles entre elles et avec \m{\sigma},
\m{\sigma'}, \m{\tau}, \m{\tau'}. Par exemple, on a
$$\sigma(z_1g_L\ot h)=\sigma(z_1\ot h)g_L, \ \ \
\sigma'(g_0h\ot z_4) = \sigma(h\ot z_4g_0), \ \ \
(g_Rt)g_L = g_R(tg_L),$$
si \ \m{z_1\in Z_1}, \m{z_4\in Z_4}, \m{h\in H}, \m{t\in T}, \m{g_L\in G_L},
\m{g_R\in G_R} et \m{g_0\in G_0}. 

\bigskip

\begin{defin}
La donn\'ee \m{\Theta} de \m{Z_1}, \m{Z_2}, \m{Z_3}, \m{Z_4}, \m{H}, \m{T},
\m{M},
\m{\sigma}, \m{\sigma'}, \m{\tau}, \m{\tau'}, et des actions de \m{G_L},
\m{G_R} et \m{G_0} s'appelle un {\em espace abstrait de complexes de type 1},
et \m{Q_C} est l'{\em espace total} de \m{\Theta}.
\end{defin}

\vskip 0.8cm

{\bf \ 3.2.1.2 }{\it Dictionnaire}

\medskip

\noindent Dans la situation du \para 3.1, on a
$$Z_1 = \Hom(\ke,\Gamma), \ \ Z_2 = \Hom(\ke,\kg),$$
$$Z_3 = \Hom(\Gamma,\kf), \ \ Z_4 = \Hom(\kg,\kf),$$
$$T = \Hom(\ke,\kf), \ \ H = \Hom(\Gamma,\kg),$$
les applications \m{\sigma}, \m{\sigma'},\m{\tau},\m{\tau'} sont les
compositions, et
$$G_L=\Aut(\ke), \ \ G_R=\Aut(\kf), \ \ G_0=\Aut(\kg).$$

\vskip 0.8cm

{\bf \ 3.2.1.3 }{\it Groupes associ\'es}

\medskip

Les groupes \m{G_L^{op}} et \m{G_R} agissent \`a gauche de mani\`ere 
\'evidente sur \m{W_C}, et \m{Q_C} est invariant par ces groupes. Soit
\m{G_1} le groupe constitu\'e des matrices
$$\left(\begin{array}{cc}g_M & 0 \cr \phi & g_0\cr\end{array}\right)$$
avec \ \m{g_M\in GL(M)}, \m{g_0\in G_0}, \m{\phi\in M^*\ot H} (la loi de
composition est \'evidente). Le groupe \m{G_1} agit lin\'eairement \`a
gauche sur \m{W_C} : cette action provient d'une action \`a gauche sur
\ \m{(Z_1\ot M)\oplus Z_2} \ et d'une action \`a droite sur
\ \m{(Z_3\ot M^*)\oplus Z_4} : si \ \m{(\phi_1,z_2,\phi_3,z_4)\in W_C} \ et
\ \m{g_1\in G_1}, on a
$$g_1(\phi_1,z_2,\phi_3,z_4) = (g_1(\phi_1,z_2), (\phi_3,z_4)g_1^{-1}).$$
L'action \`a gauche de \m{G_1} sur \ \m{(Z_1\ot M)\oplus Z_2} \ est :
$$\left(\begin{array}{cc}g_M & 0 \cr \phi & g_0\cr\end{array}\right)
\left(\begin{array}{c}\phi_1\cr z_2\end{array}\right) \  = \ 
\left(\begin{array}{cc}
g_M\phi_1\cr \sigma(\pline{\phi,\phi_1})+g_0z_2\end{array}\right).$$
L'action \`a droite de \m{G_1} sur \ \m{(Z_3\ot M^*)\oplus Z_4} \ est :
$$(\phi_3,z_4)\left(\begin{array}{cc}g_M & 0 \cr \phi & g_0\cr\end{array}
\right) \ = \ 
(\phi_3g_M+(I_M\ot \sigma')(\phi\ot z_4), z_4g_0).$$
En utilisant la commutativit\'e du diagramme \m{(D)}, on montre ais\'ement
que \m{Q_C} est \m{G_1}-invariant. Plus g\'en\'eralement, on montre que
l'application
$$W_C\lra T$$
$$(\phi_1,z_2,\phi_3,z_4)\longmapsto \tau(\pline{\phi_1,\phi_3})+
\tau'(z_2\ot z_4)$$
est \m{G_1}-invariante. Dans la situation du \para 2.2.1.2, on a
\ \m{G_1=\Aut((\Gamma\ot M)\oplus\kg)}.

On pose
$$G \ \ = \ \ G_L^{op}\times G_1\times G_R,$$
qui agit \`a gauche sur \m{W_C} et \m{Q_C}.

\vskip 1cm

\begin{subsub}{Espaces de complexes de type 2}\end{subsub}

{\bf \ 3.2.2.1 }{\it D\'efinition}

\medskip

Soient \m{Z_1}, \m{Y_2}, \m{T_2}, \m{Z_3}, \m{T}, \m{K}, \m{N} des
espaces vectoriels de dimension finie, et
$$W'_C = (Z_1\ot N)\oplus(Y_2\ot N)\oplus(Z_3\ot N^*).$$
Soient
$$\nu : K\ot Y_2\lra Z_1,$$
$$\nu' : K\ot T_2\lra T,$$
$$\lambda : Y_2\ot Z_3\lra T_2,$$
$$\tau : Z_1\ot Z_3\lra T$$
des applications lin\'eraires. On suppose que \m{\nu} induit une
inclusion \ \m{K\subset Z_1\ot Y_2^*} \ et que \m{\lambda} est surjective.
On suppose aussi que le diagramme suivant est commutatif :
$$(D') \ \ \ \ \ \ \begin{CD}
K\ot Y_2\ot Z_3  @>{\nu\ot I_{Z_3}}>> Z_1\ot Z_3\cr
@V{I_K\ot\lambda}VV @VV{\tau}V\cr
K\ot T_2 @>>{\nu'}>  T\end{CD}$$
Soit \ \m{Q'_C\subset W'_C} \ l'ensemble des points 
\m{(\psi_1,\psi_2,\psi_3)} tels que
$$\tau(\pline{\psi_1,\psi_3}) \ = \lambda(\pline{\psi_2,\psi_3}) \ = \ 0,$$
o\`u \m{\pline{\ \ }} d\'esigne la contraction de \m{N}. 

Soient \m{G_L}, \m{G_0}, \m{G_R} des groupes. On suppose que
\begin{itemize}
\item[] $G_L$ op\`ere lin\'eairement \`a droite sur $T$, $Z_1$ et $K$,
\item[] $G_R$ op\`ere lin\'eairement \`a gauche sur $T$, $Z_3$ et $T_2$,
\item[] $G_0$ op\`ere lin\'eairement \`a gauche sur $K$ et
lin\'eairement \`a droite sur $T_2$ et $Y_2$.
\end{itemize}
On suppose comme pour les complexes de type 1 que les actions des groupes
sont compatibles entre elles et avec les applications \m{\nu}, \m{\nu'},
\m{\lambda} et \m{\tau}.

\bigskip

\begin{defin}
La donn\'ee \m{\Theta'} de \m{Z_1}, \m{Y_2}, \m{T_2}, \m{Z_3}, \m{T}, \m{K},
\m{N}, \m{\nu}, \m{\nu'}, \m{\lambda}, \m{\tau} et des actions de
\m{G_L}, \m{G_0}, \m{G_R} s'appelle un {\em espace abstrait de complexes
de type 2}, et \m{Q'_C} est {\em l'espace total de } \m{\Theta'}.
\end{defin}

\vskip 0.8cm

{\bf \ 3.2.2.2 }{\it Dictionnaire}

\medskip

Dans la situation du \para 3.1, on a
$$Z_1=\Hom(\ke,\Gamma), \ Y_2=\Hom(\ke',\Gamma),$$
$$Z_3 = \Hom(\Gamma,\kf), \ \ T_2=\Hom(\ke',\kf),$$
$$T=\Hom(\ke,\kf), \ \ K=\Hom(\ke,\ke'),$$
les applications \m{\nu}, \m{\nu'}, \m{\lambda}, \m{\tau} sont les
compositions et
$$G_L=\Aut(\ke), \ G_R=\Aut(\kf), \ G_0=\Aut(\ke').$$

\vskip 0.8cm

{\bf \ 3.2.2.3 }{\it Groupes associ\'es}

\medskip

Les groupes \m{G_L^{op}} et \m{G_R} agissent \`a gauche de mani\`ere
\'evidente sur \m{W'_C}, et \m{Q'_C} est invariant par ces groupes. Soit
\m{G'_1} le groupe constitu\'e des matrices 
$$\left(\begin{array}{cc}g_L & 0\cr k & g_0\end{array}\right)$$
avec \ \m{g_L\in G_L}, \m{g_0\in G_0}, \m{k\in K} (la loi de composition
est \'evidente). Alors \m{G'_1} agit \`a droite sur \ \m{Z_1\oplus Y_2} :
$$(z_1,y_2)\left(\begin{array}{cc}g_L & 0\cr k & g_0\end{array}\right) \ = \
(z_1g_L+\nu(k\ot y_2),y_2g_0).$$
Dans la situation du \para 3.2.2.2, on a \ \m{G'_1=\Aut(\ke\oplus\ke')}.

On en d\'eduit une action \`a gauche de
$$G' \ = \ GL(N)\times {G'_1}^{op}\times G_R$$
sur \m{W'_C}. On v\'erifie comme dans le cas des complexes de type 1 que
\m{Q'_C} est \m{G'}-invariant.

\vskip 1cm

\begin{subsub}{Mutations \ 1 \m{\Longrightarrow} 2}\end{subsub}

{\bf \ 3.2.3.1 }{\it Mutations d'espaces abstraits de complexes}

\medskip

On consid\`ere l'espace abstrait de complexes de type 1 \m{\Theta} du 
\para 3.2.1. On va en d\'eduire \m{\Theta'}, espace abstrait de complexes de
type 2.
Les espaces vectoriels \m{Z_1}, \m{Z_3}, \m{T} de \m{\Theta'} sont les
m\^emes que ceux de \m{\Theta}. On prend
$$N=M\oplus H, \ \ Y_2=H^*, \ \ T_2=(H^*\ot Z_3)/Z_4, \ \ 
K=\ker(\sigma)\subset Z_1\ot H.$$
L'application \m{\tau} de \m{\Theta'} est la m\^eme que celle de \m{\Theta}.

L'application \ \m{\lambda : Y_2\ot Z_3\lra T_2} \ est la projection \
\m{H^*\ot Z_3\lra (H^*\ot Z_3)/Z_4}. 

L'application \ \m{\nu : K\ot Y_2\lra Z_1} \ est la compos\'ee
$$K\ot Y_2\subset Z_1\ot H\ot Y_2 = Z_1\ot H\ot H^*\lra Z_1.$$

Pour d\'efinir \ \m{\nu' : K\ot T_2\lra T} \ on part du
diagramme commutatif suivant, d\'eduit de \m{(D')}
$$\begin{CD}Z_1\ot H\ot Z_4 @>{\sigma\ot I_{Z_4}}>> Z_2\ot Z_4\cr
@V{\alpha}VV @VV{\tau'}V\cr
Z_1\ot H\ot Z_3\ot H^* @>>{\tau\ot tr}> T\cr\end{CD}$$
(\m{tr} d\'esignant la trace \ \m{H\ot H^*\lra\C^{}} \ et \m{\alpha} 
provenant de l'inclusion \ \m{Z_4\subset Z_3\ot H^*} \ d\'eduite de
\m{\sigma'}). Il en d\'ecoule que
$$(\tau\ot tr)\circ\alpha(\ker(\sigma)\ot Z_4) \ = \ \nsp .$$
Donc \ \m{\tau\ot tr} \ induit une application lin\'eaire
$$\ker(\sigma)\ot((Z_3\ot H^*)/Z_4) \ = \ K\ot T_2\lra T$$
qui est par d\'efinition \m{\nu'}. 

La commutativit\'e de \m{(D')} se v\'erifie ais\'ement. 

Les groupes \m{G_R}, \m{G_L}, \m{G_0} de \m{\Theta'} sont les m\^emes
que ceux de \m{\Theta} et leurs actions sont \'evidentes. L'espace
abstrait de complexes de type 2 \m{\Theta'} est ainsi compl\`etement
d\'efini. On notera
$$\Theta' \ = \ D_0(\Theta).$$

\vskip 0.8cm

{\bf \ 3.2.3.2 }{\it Mutations de complexes}

\medskip

Soit \ \m{(\phi_1,z_2,\phi_3,z_4)\in Q_C}. On va en d\'eduire une orbite
\ \m{G'.(\psi_1,\psi_2,\psi_3)} \ de \m{Q'_C}.
\begin{itemize}
\item[--] D\'efinition de \ $\psi_1\in Z_1\ot N$ : on
prend \ $\psi_0\in Z_1\ot H$ \ tel que \ $\sigma(\psi_0)=z_2$, et
$$\psi_1=\psi_0+\phi_1 \ \in \ (Z_1\ot H)\oplus(Z_1\ot M)=Z_1\ot N.$$
\item[--] D\'efinition de \ $\psi_2\in Y_2\ot N$ : on prend
$$\psi_2=I_H \ \in \ H^*\ot H\subset H^*\ot(M\oplus H)=Y_2\ot N.$$
\item[--] D\'efinition de \ $\psi_3\in Z_3\ot N^*$ : on prend
$$\psi_3=\phi_3+z_4 \ \in \ (Z_3\ot M^*)\oplus Z_4\subset
(Z_3\ot M^*)\oplus (Z_3\ot H^*)=Z_3\ot N^*.$$
\end{itemize}

\medskip

On v\'erifie ais\'ement que \m{(\psi_1,\psi_2,\psi_3)} est un \'el\'ement de
\m{Q_C} d\'efini \`a l'action pr\`es du sous-groupe de \m{G'_1} isomorphe 
\`a \m{K}, constitu\'e des matrices
$$\left(\begin{array}{cc}1 & 0\cr k & 1\end{array}\right), \ \ \ k\in K.$$

On a donc d\'efini une application
$$\ov{D_0} : Q_C\lra Q'_C/G' .$$

\bigskip

\begin{xlemm}
Si \ \m{x\in Q_C} \ et \ \m{g\in G}, alors \ \m{\ov{D_0}(gx)=\ov{D_0}(x)}.
Donc \m{\ov{D_0}} induit une application
$$D_0 : Q_C/G\lra Q'_C/G'.$$
\end{xlemm}

\medskip

V\'erification imm\'ediate. 

\vskip 1cm

\begin{subsub}{Mutations \ 2 \m{\Longrightarrow} 1}\end{subsub}

{\bf \ 3.2.4.1 }{\it Mutations d'espaces abstraits de complexes}

\medskip

On consid\`ere l'espace abstrait de complexes de type 2 \m{\Theta'} du
\para 3.2.2. On va en d\'eduire \m{\Theta}, espace abstrait de complexes
de type 1. On doit supposer que
\ \m{\dim(N)\geq\dim(Y_2)}.
On note \m{{W'}^0_C} le sous-ensemble \m{G'}-invariant de \m{W'_C}
constitu\'e des \m{(\psi_1,\psi_2,\psi_3)} tels que \ 
\m{\psi_2 : Y_2^*\lra N} \ soit injective. Soit \ \m{{Q'}^0_C=
Q'_C\cap {W'}^0_C}.

On d\'efinit maintenant \m{\Theta}. Les espaces vectoriels \m{Z_1}, \m{Z_3},
\m{T} de \m{\Theta} sont les m\^emes que ceux de \m{\Theta'}. On prend pour
\m{M} un espace vectoriel de dimension \ \m{\dim(N)-\dim(Y_2)}, et
$$H=Y_2^*, \ \ Z_2=(Z_1\ot Y_2^*)/K, \ \ Z_4=\ker(\lambda)\subset 
Y_2\ot Z_3.$$

\noindent L'application \m{\tau} de \m{\Theta} est la m\^eme que celle de
\m{\Theta'}.

\noindent L'application \ \m{\sigma : Z_1\ot H\lra Z_2} \ est la projection
\ \m{Z_1\ot Y_2^*\lra (Z_1\ot Y_2^*)/K}.

\noindent L'application \ \m{\sigma' : H\ot Z_4\lra Z_3} \ est la restriction
de
$$tr\ot I_{Z_3} : H\ot H^*\ot Z_3\lra Z_3.$$

Pour d\'efinir \ \m{\tau' : Z_2\ot Z_4\lra T} \ on part du 
diagramme commutatif suivant d\'eduit de \m{(D')}
$$\begin{CD}K\ot Y_2\ot Z_3 @>{I_K\ot\lambda}>> K\ot T_2\cr
@V{\ov{\nu}\ot I_{Y_2\ot Z_3}}VV @VV{\nu'}V\cr
Z_1\ot Y_2^*\ot Y_2\ot Z_3 @>>{\tau\ot tr}> T\cr\end{CD}$$
\m{\ov{\nu}} d\'esignant l'inclusion \ \m{K\subset Z_1\ot Y_2^*} \ d\'eduite
de \m{\nu}. On a donc
$$(\tau\ot tr)\circ(\ov{\nu}\ot I_{Y_2\ot Z_3})(K\ot \ker(\lambda)=\nsp,$$
et \ \m{\tau\ot tr} \ induit donc
$$\tau' : ((Z_1\ot Y_2^*)/K)\ot\ker(\lambda)=Z_2\ot Z_4\lra T.$$
La commutativit\'e du diagramme \m{(D)} se v\'erifie ais\'ement. 
Les groupes \m{G_0}, \m{G_L} et \m{G_R} de \m{\Theta} sont les m\^emes que
ceux de \m{\Theta'} et leurs actions sont \'evidentes. L'espace abstrait de
complexes de type 1 \m{\Theta} est ainsi compl\`etement d\'efini. On notera
$$\Theta \ = \ D'_0(\Theta').$$

\bigskip

\begin{xprop}
On a \ \m{D_0\circ D'_0(\Theta')=\Theta'} \ et \ 
\m{D'_0\circ D_0(\Theta)=\Theta}.
\end{xprop}

\medskip

Imm\'ediat. 

\vskip 0.8cm

{\bf \ 3.2.4.2 }{\it Mutations de complexes}

\medskip

Soit \ \m{(\psi_1,\psi_2,\psi_3)\in {Q'}^0_C}. On va en d\'eduire une
orbite \ \m{G.(\phi_1,z_2,\phi_3,z_4)} \ de \m{Q_C}. 
\begin{itemize}
\item[--] D\'efinition de \ $\phi_1\in Z_1\ot M$ : on fixe d'abord un
isomorphisme entre $M$ et un suppl\'ementaire de l'image de $\psi_2$ dans
$N$ : 
$$N \ = \ Y_2^*\oplus M.$$
On prend pour $\phi_1$ la composante de $\psi_1$ dans \ $Z_1\ot M$.
\item[--] D\'efinition de \ $z_2\in Z_2$ : on prend la projection sur \
$(Z_1\ot Y_2^*)/K$ \ de la composante de $\psi_2$ dans \ $Z_1\ot Y_2^*$.
\item[--] D\'efinition de \ $\phi_3\in Z_3\in M^*$ : on prend la composante
de $\psi_3$ dans \ $Z_3\ot M^*$.
\item[--] D\'efinition de \ $z_4\in Z_4$ : on prend \
$z_4=\pline{\psi_2,\psi_3}\in\ker(\lambda).$
\end{itemize}

\medskip

On v\'erifie ais\'ement que \m{(\phi_1,z_2,\phi_3,z_4)} est un \'el\'ement 
de \m{Q_C} d\'efini \`a l'action pr\`es du groupe \m{G_1}.

On a donc d\'efini une application
$$\ov{D'_0} : {Q'}^0_C\lra Q_C/G .$$

\bigskip

\begin{xlemm}
Si \ \m{x\in {Q'}^0_C} \ et \ \m{g\in G'}, alors
 \ \m{\ov{D'_0}(gx)=\ov{D'_0}(x)}.
Donc \m{\ov{D'_0}} induit une application
$$D'_0 : {Q'}^0_C/G'\lra Q_C/G.$$
\end{xlemm}

\medskip

V\'erification imm\'ediate. 

\vskip 1cm

\begin{subsub}{Th\'eor\`emes d'isomorphisme}\end{subsub}

\begin{xtheo}
On a \ \m{D'_0\circ D_0 = I_{Q_C/G}} \ et \ 
\m{D_0\circ D'_0 = I_{{Q'}^0_C/G'}}.
\end{xtheo}

\medskip

V\'erification imm\'ediate. 

\bigskip

On a donc obtenu une bijection canonique
$$Q_C/G \ \simeq {Q'}^0_C/G'.$$

\medskip

On suppose maintenant que les groupes sont alg\'ebriques, ainsi que leurs
actions sur les espaces vectoriels dont il est question. On a alors :

\bigskip

\begin{xtheo}
Il existe un quasi-isomorphisme fort canonique de \m{Q_C} vers \m{{Q'}^0_C},
au dessus de l'isomorphisme pr\'ec\'edent.
\end{xtheo}

\bigskip

\begin{proof}
Le quasi-isomorphisme fort \ \m{Q_C\lra {Q'}^0_C} \ est d\'efini par une
seule carte, obtenue en utilisant une section de \m{\sigma}. Le 
quasi-isomorphisme inverse est d\'efini par des cartes ind\'ex\'ees sur
la grassmannienne \m{Gr_0} des sous-espaces vectoriels de $N$ de dimension
\'egale \`a celle de $M$. L'ouvert correspondant \`a \ \m{M_0\in Gr_0} \
est l'ensemble des points \m{(\psi_1,\psi_2,\psi_3)} de \m{{Q'}^0_C} tels que
l'image de \ \m{\psi_2 : Y_2^*\lra N} \ ne rencontre pas \m{M_0}. Les
v\'erifications (fastidieuses) sont laiss\'ees au lecteur.
\end{proof}

\bigskip

On peut donner une version abstraite de la proposition 3.3, et obtenir des
r\'esultats analogues aux th\'eor\`emes 3.7 et 3.8.

\vskip 1.5cm

\begin{sub}{\bf Application aux espaces de complexes}\end{sub}

Soient $X$
une vari\'et\'e alg\'ebrique projective, \m{p\geq 1} \ un entier, 
\m{n_0,\ldots,n_p} des entiers positifs, et pour \ \m{0\leq i\leq p}, 
\m{1\leq j\leq n_i}, \ \m{\ke^{(i)}_j} un faisceau coh\'erent sur $X$ et 
\m{M^{(i)}_j} un espace vectoriel non nul de dimension finie. On pose, pour
\m{0\leq i\leq p}
$$\ke_i \ = \ \som_{1\leq j\leq n_i}(\ke^{(i)}_j\ot M^{(i)}_j).$$
On suppose que les faisceaux \m{\ke^{(i)}_j} sont simples, et que
$$\Hom(\ke^{(i)}_j, \ke^{(i')}_{j'}) \ = \ \nsp$$
si \ \m{i>i'}, ou \ \m{i=i'}, \m{j>j'}. Soit \m{Q_C} la vari\'et\'e des
complexes
$$\ke_0\lra\ke_1\lra\ldots\lra\ke_p.$$
On pose \ \m{\ke_{i}=0} \ si \m{i<0} ou \m{i>p}.
Soit \m{i_0} un entier, avec \ \m{0\leq i_0\leq p}. On pose
$$\ke=\ke_{i_0-1}, \ \ \Gamma=\ke^{(i_0)}_1, \ \ M=M^{(i_0)}_1, \ \
\kg=\som_{2\leq j\leq n_i}(\ke^{(i_0)}_j\ot M^{(i_0)}_j), \ \
\kf=\ke_{i_0+1},$$
de telle sorte que les complexes pr\'ec\'edents se mettent sous la forme
$$\cdots\lra\ke_{i_0-2}\lra\ke\lra(\Gamma\ot M)\oplus\kg\lra\kf\lra\ke_{i_0+2}\lra
\cdots$$
On suppose que les conditions du \para 3.1.1 sont v\'erifi\'ees, ce qui 
d\'efinit le faisceau \m{\ke'}. Soit \m{Q'_C} la vari\'et\'e des 
complexes du type
$$\cdots\lra\ke_{i_0-2}\lra\ke\oplus\ke'\lra\Gamma\ot N\lra\kf\lra\ke_{i_0+2}\lra
\cdots$$
et \m{{Q'}^0_C} l'ouvert de \m{Q'_C} constitu\'e des complexes tels que le
morphisme \ \m{\ke'\lra\Gamma\ot N} \ induise une inclusion
\ \m{\Hom(\ke',\Gamma)^*\subset N}. 
Le groupe \m{G} (resp. \m{G'}) op\'erant sur \m{Q_C} (resp. \m{Q'_C}) est
le produit des groupes d'automorphismes des termes des complexes. On 
g\'en\'eralise sans difficult\'es les r\'esultats du \para 3.2 \`a ce cas et
on d\'efinit de mani\`ere \'evidente les mutations de complexes de \m{Q_C} ou
\m{{Q'}^0_C} (qui sont des complexes de \m{{Q'}^0_C} et \m{Q_C} respectivement),
et on obtient l'analogue du th\'eor\`eme 3.8 pour les vari\'et\'es de
complexes.

\vskip 2.5cm

\section{Vari\'et\'es de modules de complexes}

La construction des vari\'et\'es de modules de morphismes de \cite{dr_tr}
est une application du th\'eor\`eme 3.8. On donne ici une autre
application de ce th\'eor\`eme \`a la construction de vari\'et\'es de
modules de complexes. On travaille dans le language des faisceaux, ce 
qui donne des d\'emonstrations plus explicites. Il est \'evidemment
possible de faire une version abstraite de la construction des 
vari\'et\'es de modules, comme dans \cite{dr_tr}.

Soient \m{\ke_1}, \m{\kf_1}, \m{\kf_2}, \m{\kg_1} des faisceaux coh\'erents 
sur une vari\'et\'e projective \m{X}, et \m{L_1}, \m{M_1}, \m{M_2}, \m{N_1}
des espaces vectoriels de dimension finie. On s'int\'eresse \`a des
complexes du type
$$\ke_1\ot L_1\lra (\kf_1\ot M_1)\oplus(\kf_2\ot M_2)\lra \kg_1\ot N_1.$$
On fait les hypoth\`eses suivantes :

\begin{itemize}
\item[--] Les faisceaux $\ke_1$, $\kf_1$, $\kf_2$ et $\kg_1$ sont simples, et
\ $\Hom(\kf_2,\kf_1)=\Hom(\kf_1,\ke_1)=\nsp.$
\item[--] Le morphisme canonique
$$\kf_1\ot\Hom(\kf_1,\kf_2)\lra\kf_2$$
est surjectif. On note \m{\kh_1} son noyau.
\item[--] Le morphisme canonique
$$\ke_1\ot\Hom(\ke_1,\kh_1)\lra\kh_1$$
est surjectif. On note \m{\kk_1} son noyau.
\item[--] On a 
$$\Ext^1(\kf_2,\ke_1)=\Ext^1(\ke_1,\kh_1)=\Ext^1(\kf_1,\kf_2)=\nsp,$$
$$\Ext^1(\kf_2,\kf_1)=\Ext^1(\kh_1,\kh_1)=\Ext^1(\kh_1,\ke_1)=\nsp.$$
\end{itemize}

Cette derni\`ere hypoth\`ese entra\^ine que les faisceaux \m{\kh_1} et 
\m{\kk_1} sont simples.

On se trouve dans la situation
de la proposition 3.1 (avec \ \m{\ke=\ke_1\ot L_1},
\m{\Gamma=\kf_1}, \m{M=M_1}, \m{\kg=\kf_2\ot M_2} et \m{\kf=\kg_1\ot N_1}). On
note \m{Q_C} la vari\'et\'e des complexes du type pr\'ec\'edent. Elle peut
donc \^etre vue comme l'espace total d'un espace abstrait de complexes de 
type 1. Le groupe \m{G} op\'erant sur \m{Q_C} est
$$G \ = \ GL(L_1)\times \Aut((\kf_1\ot M_1)\oplus(\kf_2\ot M_2))\times
GL(N_1).$$
On peut voir \ \m{\Aut((\kf_1\ot M_1)\oplus(\kf_2\ot M_2))} \ comme constitu\'e
de matrices du type
$$\left(\begin{array}{cc}g_1 & 0\cr\phi & g_2\cr\end{array}\right),$$
avec  \ \m{g_1\in GL(M_1)}, \m{g_2\in GL(M_2)} \ et \ \m{\phi\in
\Hom(\Hom(\kf_1,\kf_2)^*\ot M_1,M_2)}.  

On pose
$$a=\dim(\Hom(\kf_1,\kf_2)),\ \ \ b=\dim(\Hom(\ke_1,\kh_1)).$$

\vskip 1.5cm

\begin{sub}{\bf Notions de (semi-)stabilit\'e et vari\'et\'es de modules de
complexes}\end{sub}

Le groupe $G$ poss\`ede deux sous-groupes importants : le premier est le
sous-groupe normal unipotent maximal \'evident $H$, isomorphe au groupe
additif 
\m{\Hom(\Hom(\kf_1,\kf_2)^*\ot M_1,M_2))}. Le second est le sous-groupe
r\'eductif
$$G_{red} \ = \ GL(L_1)\times GL(M_1)\times GL(M_2)\times GL(N_1),$$
dont l'inclusion dans $G$ induit un isomorphisme \ 
\m{G_{red}\simeq G/H}.

L'action de \m{G_{red}} est un cas particulier des actions \'etudi\'ees 
dans \cite{king}. On consid\`ere l'espace vectoriel
$$W \ = \ \Hom(\ke_1\ot L_1,(\kf_1\ot M_1)\oplus(\kf_2\ot M_2))\times
\Hom((\kf_1\ot M_1)\oplus(\kf_2\ot M_2),\kg_1\ot N_1),$$
dont $Q_C$ est une sous-vari\'et\'e ferm\'ee. L'action de $G$ sur \m{Q_C}
s\'etend de mani\`ere \'evidente \`a $W$.

Soient \m{\lambda_1}, \m{\mu_1}, \m{\mu_2}, \m{\nu_1} des nombres rationnels
non nuls tels que
$$\lambda_1\dim(L_1)+\mu_1\dim(M_1)+\mu_2\dim(M_2)+\nu_1\dim(N_1)=0.$$ 

\bigskip

\begin{defin}
Un point \m{(\phi,\psi)} de $W$ est dit {\em $G_{red}$-semi-stable} (resp.
{\em $G_{red}$-stable}) relativement \`a \m{(\lambda_1,\mu_1,\mu_2,\nu_1)} si
pour tous sous-espaces vectoriels \ \m{L'_1\subset L_1}, \m{M'_1\subset M_1},
\m{M'_2\subset M_2}, \m{N'_1\subset N_1}, avec \ 
\m{(L'_1,M'_1,M'_2,N'_1)\not = (L_1,M_1,M_2,N_1)} \ ou 
\m{(\nsp,\nsp,\nsp,\nsp)}, tels que
$$\phi(\ke_1\ot L'_1)\subset (\kf_1\ot M'_1)\oplus(\kf_2\ot M'_2), \ \ \
\psi((\kf_1\ot M'_1)\oplus(\kf_2\ot M'_2))\subset \kg_1\ot N'_1,$$
on a
$$\lambda_1\dim(L'_1)+\mu_1\dim(M'_1)+\mu_2\dim(M'_2)+\nu_1\dim(N'_1)\leq 0
\ \ \ {\rm (resp. \ \ } \ < \ {\rm )}.$$
\end{defin}

\bigskip 

On dit que \m{(\lambda_1, \mu_1, \mu_2, \nu_1)} est une {\em polarisation}
de l'action de $G$ sur $W$ (ou \m{Q_C}).
On note \m{W^{ss}_{red}} (resp. \m{W^s_{red}}) l'ouvert de \m{W} constitu\'e
des points \m{G_{red}}-semi-stables (resp. \m{G_{red}}-stables). Soient
\ \m{Q^{ss}_{C,red}=W^{ss}_{red}\cap Q_C}, 
\m{Q^s_{C,red}=W^s_{red}\cap Q_C}. D'apr\`es \cite{king}, il existe un
bon quotient \ \m{W^{ss}_{red}//G} \ et un quotient g\'eom\'etrique
lisse \m{W^s_{red}/G}. Par cons\'equent il existe aussi un bon quotient
\ \m{Q^{ss}_{C,red}//G}. Mais ce n'est pas le quotient que nous
recherchons.

\bigskip

\begin{defin}
Un point \m{x} de $W$ est dit {\em $G$-semi-stable} (resp.
{\em $G$-stable}) relativement \`a \m{(\lambda_1,\mu_1,\mu_2,\nu_1)} si
tous les points de l'orbite \m{H.x} sont $G_{red}$-semi-stables (resp.
$G_{red}$-stables).
\end{defin}

\bigskip 

On note \m{W^{ss}} (resp. \m{W^s}) l'ouvert de \m{W} constitu\'e
des points \m{G}-semi-stables (resp. \m{G}-stables). Soient
\ \m{Q^{ss}_{C}=W^{ss}\cap Q_C}, \m{Q^s_{C}=W^s\cap Q_C}. On cherche
\`a prouver l'existence de quasi-bons quotients \ \m{Q^{ss}_C//G}. De tels 
quotients seront appel\'es des {\em vari\'et\'es de modules de complexes}.

\medskip

On montre ais\'ement que si \m{W^{s}} est non vide on doit avoir
$$\lambda_1>0, \ \ \nu_1<0, \ \ \mu_1\dim(M_1)+\nu_1\dim(N_1)<0, \ \
\mu_2\dim(M_2)+\nu_1\dim(N_1)<0.$$
On supposera par la suite que ces in\'egalit\'es sont v\'erifi\'ees.

\vskip 1.5cm

\begin{sub}{\bf Mutations et polarisations associ\'ees}\end{sub}

\begin{subsub}{La premi\`ere mutation}\end{subsub}

{\bf \ 4.2.1.1 - }{\it D\'efinition}

\medskip

Soit
$$(*) \ \ \ \ \ \ \ \ke_1\ot L_1 \ \hfl{(f_1,f_2)}{} \  
(\kf_1\ot M_1)\oplus(\kf_2\ot M_2) \ \hfl{(g_1,g_2)}{} \ \kg_1\ot N_1.$$
un complexe. Une premi\`ere mutation donne un complexe
$$(**) \ \ \ \ \ \ \ (\ke_1\ot L_1)\oplus(\kh_1\ot M_2) \ 
\hfl{(\phi_1,\phi_2)}{} \ 
\kf_1\ot P_1\ \hfl{\phi}{} \ \kg_1\ot N_1,$$
avec
$$P_1=(\Hom(\kf_1,\kf_2)\ot M_2)\oplus M_1.$$
Le groupe op\'erant sur la vari\'et\'e \m{Q'_C} de ces complexes est
$$G'=\Aut((\ke_1\ot L_1)\oplus(\kh_1\ot M_2))\times GL(P_1)\times GL(N_1).$$
Le groupe \ \m{\Aut((\ke_1\ot L_1)\oplus(\kh_1\ot M_2))} \ est constitu\'e de
matrices
$$\left(\begin{array}{cc}g_1 & 0\cr \phi & g_2\cr\end{array}\right),$$
avec \ \m{g_1\in GL(L_1)}, \m{g_2\in GL(M_2)}, 
\m{\phi\in\Hom(\Hom(\ke_1,\kh_1)^*\ot L_1,M_2)} (on a 
\m{\Hom(\kh_1,\ke_1)=\nsp}, \`a cause du fait que
\m{\Ext^1(\kf_2,\ke_1)=\Hom(\kf_1,\ke_1)=\nsp}). Le complexe \m{(**)} n'est
pas unique, mais sa $G'$-orbite l'est.

D\'ecrivons maintenant une mutation du complexe \m{(*)}. Le morphisme 
$$\phi_1 : \ke_1\ot L_1\lra \kf_1\ot P_1$$
est la somme de 
$$f_1 : \ke_1\ot L_1\lra \kf_1\ot M_1$$
et d'un rel\`evement de \ \m{f_2 : \ke_1\ot L_1\lra \kf_2\ot M_2} \ en un
morphisme \hfil\break
\m{\ke_1\ot L_1\lra \kf_1\ot\Hom(\kf_1,\kf_2)\ot M_2}. Un tel 
rel\`evement est possible car 
\ \m{\Ext^1(\ke_1,\kh_1)=\nsp}. Le morphisme
$$\phi_2 : \kh_1\ot M_2\lra\kf_1\ot P_1$$
est \'egal \`a \ \m{\sigma\ot I_{M_2}}, o\`u \m{\sigma} est l'inclusion \
\m{\kh_1\lra\kf_1\ot\Hom(\kf_1,\kf_2)}. Le morphisme \m{\phi} est \'egal \`a 
\m{g_1} sur \m{\kf_1\ot M_1}, et sur \m{\kf_1\ot \Hom(\kf_1,\kf_2)\ot M_2}, c'est
la compos\'ee
$$\kf_1\ot \Hom(\kf_1,\kf_2)\ot M_2\lra \kf_2\ot M_2 \ \hfl{g_2}{} \ \kg_1\ot N_1.$$

\vskip 0.8cm

{\bf \ 4.2.1.2 - }{\it (Semi-)stabilit\'e dans \m{Q'_C}}

\medskip

On d\'efinit comme pour les complexes de type \m{(*)} une notion de 
(semi-)stabilit\'e, d\'ependant d'une suite \m{(\alpha_1, \alpha_2,\beta_1,
\gamma_1)} de nombres rationnels telle que
$$\alpha_1\dim(L_1)+\alpha_2\dim(M_2)+\beta_1\dim(P_1)+\gamma_1\dim(N_1)
=0.$$ 
Soient \m{G'_{red}} les sous-groupe r\'eductif canonique de $G'$,
et $H'$ le sous-groupe normal unipotent maximal isomorphe au groupe additif
\m{\Hom(L_1\ot\Hom(\ke_1,\kh_1),M_2)}.
La \m{G'_{red}}-\hbox{(semi-)}stabilit\'e relativement \`a 
\m{(\alpha_1,\alpha_2,\beta_1,\gamma_1)} 
est encore un cas particulier des actions 
\'etudi\'ees dans \cite{king} : le complexe \m{(**)} est 
\m{G'_{red}}-semi-stable (resp. \m{G'_{red}}-stable) si pour tous 
sous-espaces vectoriels \ \m{L'_1\subset L_1}, \m{M'_2\subset M_2},
\m{P'_1\subset P_1}, \m{N'_1\subset N_1}, avec 
\ \m{(L'_1,M'_2,P'_1,N'_1)\not = (L_1,M_2,P_1,N_1)} \ ou
\m{(\nsp,\nsp,\nsp,\nsp)}, tels que
$$(\phi_1,\phi_2)((\ke_1\ot L'_1)\oplus(\kh_1\ot M'_2))\subset\kf_1\ot P'_1,
\ \ \ \phi(\kf_1\ot P'_1)\subset\kg_1\ot N'_1,$$
on a
$$\alpha_1\dim(M'_1)+\alpha_2\dim(M'_2)+\beta_1\dim(P'_1)+
\gamma_1\dim(N'_1)\leq 0 \ \ \ \ {\rm(resp. \ } \ < \ {\rm )}.$$
Le complexe \m{(**)} est $G'$-semi-stable (resp. \m{G'}-stable) si tous les
points de sa $H'$-orbite sont \m{G'_{red}}-semi-stables (resp.
\m{G'_{red}}-stables). 

Soient \ \m{L'_1\subset L_1}, \m{M'_1\subset M_1}, \m{M'_2\subset M_2}, 
\m{N'_1\subset N_1} \ des sous-espaces vectoriels tels que
$$(f_1,f_2)(\ke_1\ot L'_1)\subset (\kf_1\ot M'_1)\oplus(\kf_2\ot M'_2), \ \ \ 
(g_1,g_2)((\kf_1\ot M'_1)\oplus(\kf_2\ot M'_2))\subset \kg_1\ot N'_1.$$
On pose 
$$P'_1=(\Hom(\kf_1,\kf_2)\ot M'_2)\oplus M'_1.$$
Le rel\`evement de \m{f_2} peut \^etre choisi de telle sorte 
qu'il envoie \m{\ke_1\ot L'_1} dans \hfil\break 
\m{\kf_1\ot\Hom(\kf_1,\kf_2)\ot M'_2}. On a alors
$$(\phi_1,\phi_2)((\ke_1\ot L'_1)\oplus(\kh_1\ot M'_2))\subset\kf_1\ot P'_1,
 \ \ \ \phi(\kf_1\ot P'_1)\subset\kg_1\ot N'_1.$$

Soit \m{(\lambda_1,\mu_1,\mu_2,\nu_1)} une polarisation de l'action de $G$
sur \m{Q_C}. On pose
$$\alpha_1=\lambda_1, \ \ \alpha_2=\mu_2-a\mu_1, \ \ \beta_1=\mu_1,
\ \ \gamma_1=\nu_1.$$
On d\'eduit imm\'ediatement de ce qui pr\'ec\`ede la

\bigskip

\begin{xprop}
Si le complexe \m{(**)} est \m{G'_{red}}-semi-stable
(resp. \m{G'_{red}}-stable) relativement \`a \m{(\alpha_1,\alpha_2,\beta_1,
\gamma_1)}, alors \m{(*)} est \m{G}-semi-stable (resp. \m{G}-stable)
relativement \`a \m{(\lambda_1,\mu_1,\mu_2,\nu_1)}.
\end{xprop}

\bigskip

On voit ais\'ement que si \m{{Q'}^s_C} est non vide, alors on a
\ \m{\mu_2 \ > \ a\mu_1}.

\vskip 1cm

\begin{subsub}{La seconde mutation}\end{subsub}

\medskip

{\bf \ 4.2.2.1 - }{\it D\'efinition}

\medskip

Une seconde mutation donne, partant du complexe \m{(**)}, 
un complexe du type 
$$(***) \ \ \ \ \ \ \ \kk_1\ot M_2 \ \hfl{\psi}{} \ \ke_1\ot Q_1
\ \hfl{\psi'}{} \ \kf_1\ot P_1 \ \hfl{\phi}{} \ \kg_1\ot N_1,$$
avec
$$Q_1=(\Hom(\ke_1,\kh_1)\ot M_2)\oplus L_1.$$
Le groupe op\'erant sur la vari\'et\'e \m{Q''_C} de ces complexes est
$$G''=GL(M_2)\times GL(Q_1)\times GL(P_1)\times GL(N_1),$$
qui est r\'eductif. Notons qu'ici la
mutation est uniquement d\'etermin\'ee (\`a partir de \m{(**)}).

D\'ecrivons maintenant une mutation du complexe \m{(**)}. 
Le morphisme \m{\psi}
est \'egal \`a \m{\sigma'\ot I_{M_2}}, o\`u \m{\sigma'} est l'inclusion \ 
\m{\kk_1\subset\ke_1\ot\Hom(\ke_1,\kh_1)}. Posons
$$\psi' \ = \ \left(\begin{array}{cc}\psi'_{11} & \psi'_{12}\cr
\psi'_{21} & \psi'_{22}\cr\end{array}\right),$$
relativement aux d\'ecompositions
$$P_1=(\Hom(\kf_1,\kf_2)\ot M_2)\oplus M_1, \ \ \
Q_1=(\Hom(\ke_1,\kh_1)\ot M_2)\oplus L_1.$$
On a 
$$\left(\begin{array}{c}\psi'_{12}\cr \psi'_{22}\cr\end{array}\right)
 \ = \ \phi_1.$$
Le morphisme \m{\psi'_{11}} provient de l'application canonique
$$\Hom(\ke_1,\kh_1)\ot\Hom(\ke_1,\kf_1)^*\lra\Hom(\kf_1,\kf_2),$$
compte tenu de l'isomorphisme canonique \ \m{\Hom(\kf_1,\kf_2)\simeq
\Hom(\kh_1,\kf_1)^*}. On a enfin \ \m{\psi'_{21}=0}. 

\vskip 0.8cm

{\bf \ 4.2.2.2 - }{\it (Semi-)stabilit\'e dans \m{Q''_C}}

\medskip

On d\'efinit comme pour les complexes de type \m{(*)} une notion de 
(semi-)stabilit\'e, d\'ependant d'une suite \m{(\delta,\epsilon,\theta,
\rho)} de nombres rationnels telle que
$$\delta\dim(M_2)+\epsilon\dim(Q_1)+\theta\dim(P_1)+\rho\dim(N_1)=0.$$ 
La \m{G''}-(semi-)stabilit\'e relativement \`a \m{(\delta, 
\epsilon,\theta,\rho)} est encore un cas particulier des actions 
\'etudi\'ees dans \cite{king} : le complexe \m{(***)} est 
\m{G''}-semi-stable (resp. \m{G''}-stable) si pour tous 
sous-espaces vectoriels \ \m{M'_2\subset M_2}, \m{Q'_1\subset Q_1},
\m{P'_1\subset P_1}, \m{N'_1\subset N_1}, avec
\m{(M'_2,Q'_1,P'_1,N'_1)\not = (M_2,Q_1,P_1,N_1)} \ ou
\m{(\nsp,\nsp,\nsp,\nsp)}, tels que
$$\psi(\kk_1\ot M'_2)\subset\ke_1\ot Q'_1, \ \psi'(\ke_1\ot Q'_1)\subset
\kf_1\ot P'_1, \ \phi(\kf_1\ot P'_1)\subset\kg_1\ot N'_1,$$
on a
$$\delta\dim(M'_2)+\epsilon\dim(Q'_1)+\theta\dim(P'_1)+
\rho\dim(N'_1)\leq 0 \ \ \ \ {\rm(resp. \ } \ < \ {\rm )}.$$

Soient \ \m{L'_1\subset L_1}, \m{M'_1\subset M_1}, \m{M'_2\subset M_2}, 
\m{N'_1\subset N_1} \ des sous-espaces vectoriels tels que
$$(f_1,f_2)(\ke_1\ot L'_1)\subset (\kf_1\ot M'_1)\oplus(\kf_2\ot M'_2), \ \ \ 
(g_1,g_2)((\kf_1\ot M'_1)\oplus(\kf_2\ot M'_2))\subset \kg_1\ot N'_1.$$
On pose 
$$P'_1=(\Hom(\kf_1,\kf_2)\ot M'_2)\oplus M'_1, \ \ 
Q'_1=(\Hom(\ke_1,\kh_1)\ot M'_2)\oplus L'_1.$$
Comme dans le \para 3.2.1.2, le rel\`evement de \m{f_2} peut \^etre choisi
de telle sorte qu'il envoie \m{\ke_1\ot L'_1} dans 
\m{\kf_1\ot\Hom(\kf_1,\kf_2)\ot M'_2}. On a alors
$$\psi(\kk_1\ot M'_2)\subset\ke_1\ot Q'_1, \ \psi'(\ke_1\ot Q'_1)\subset
\kf_1\ot P'_1, \ \phi(\kf_1\ot P'_1)\subset\kg_1\ot N'_1.$$

Soit \m{(\lambda_1,\mu_1,\mu_2,\nu_1)} une polarisation de l'action de $G$
sur \m{Q_C}. On pose
$$\delta=\mu_2-a\mu_1-b\lambda_1, \ \ \epsilon=\lambda_1, \ \ \theta=\mu_1,
\ \ \rho=\nu_1.$$
On d\'eduit imm\'ediatement de ce qui pr\'ec\`ede la

\bigskip

\begin{xprop}
Si le complexe \m{(***)} est \m{G''}-semi-stable
(resp. \m{G''}-stable) relativement \`a \m{(\delta,\epsilon,\theta,
\rho)}, alors \m{(*)} est \m{G}-semi-stable (resp. \m{G}-stable)
relativement \`a \m{(\lambda_1,\mu_1,\mu_2,\nu_1)}.
\end{xprop}

On voit ais\'ement que si \m{{Q''}^s_C} est non vide, alors on a
\ \m{\mu_2 \ > \ a\mu_1+b\lambda_1}.
Notons que si \m{Q_C^s} est non vide, cette condition est plus forte 
que celle que l'on avait trouv\'ee en supposant que \m{{Q'}^s_C} est non
vide.

\vskip 1.5cm

\begin{sub}{\bf Cas d'\'equivalence des (semi-)stabilit\'es et construction des
vari\'et\'es de modules}\end{sub}

\begin{subsub}{D\'efinitions de constantes}\end{subsub}

On d\'efinit ici des constantes qui interviendront par la suite. Elles
ont d\'ej\`a \'et\'e d\'efinies et utilis\'ees sans \cite{dr_tr}.
On consid\`ere l'application canonique
$$\tau : \Hom(\kf_1,\kg_1)^*\ot\Hom(\kf_1,\kf_2)\lra\Hom(\kf_2,\kg_1)^*.$$
Pour tout entier positif $k$, soit
$$\tau_k = \tau\ot I_{\C^k} :
\Hom(\kf_1,\kg_1)^*\ot(\Hom(\kf_1,\kf_2)\ot\C^{k})
\lra\Hom(\kf_2,\kg_1)^*\ot\C^{k}.$$
Soit \m{G_k} l'ensemble des sous-espaces vectoriels propres \m{K} de 
\m{\Hom(\kf_1,\kf_2)\ot\C^{k}} tels que pour tout sous-espace vectoriel
propre \ \m{V\subset\C^{k}}, \m{K} ne soit pas contenu dans
\m{\Hom(\kf_1,\kf_2)\ot V}. On pose
$$c_1(k)= \supp_{K\in{G_k}}(\frac{\codim(
\tau_k(\Hom(\kf_1,\kf_2)\ot K))}{\codim(K)}).$$
On d\'efinit de m\^eme \m{c_2(k)}, qui correspond \`a l'application
canonique
$$\tau' : \Hom(\ke_1,\kf_1)^*\ot\Hom(\ke_1,\kh_1)\lra\Hom(\kf_1,\kf_2),$$
provenant de l'isomorphisme canonique \
\m{\Hom(\kf_1,\kf_2)\simeq\Hom(\kh_1,\kf_1)^*}.

Il est clair qu'on a
\ \m{c_1(k+1)\geq c_1(k), \ \ c_2(k+1)\geq c_2(k)}.

\vskip 1cm

\begin{subsub}{Equivalence des (semi-)stabilit\'es}\end{subsub}

\begin{xprop}
On suppose que
$$\mu_2 \ \geq \ (b-\frac{a}{c_2(m_2)})\lambda_1-ac_1(m_2)\nu_1.$$
Alors le complexe \m{(*)} est \m{G}-semi-stable relativement \`a
\m{(\lambda_1,\mu_1,\mu_2,\nu_1)} si et seulement si \m{(***)} est
\m{G''}-semi-stable relativement \`a \m{(\delta,\epsilon,\theta,\rho)}.
\end{xprop}

\begin{proof} D'apr\`es la proposition 4.2, il suffit de montrer que si \m{(***)}
n'est pas \m{G''}-semi-stable, \m{(*)} n'est pas \m{G}-semi-stable. 

On suppose que le complexe \m{(***)} n'est pas semi-stable. 
Soient \ \m{M'_2\subset M_2}, \m{Q'_1\subset Q_1}, \m{P'_1\subset P_1},
\m{N'_1\subset N_1} \ des sous-espaces vectoriels, de dimensions respectives 
\m{m'_2}, \m{q'_1}, \m{p'_1}, \m{n'_1}, tels que
$$\psi(\kk_1\ot M'_2)\subset\ke_1\ot Q'_1, \ \
\psi'(\ke_1\ot Q'_1)\subset\kf_1\ot P'_1, \ \ \phi(\kf_1\ot P'_1)\subset
\kg_1\ot N'_1,$$
et
$$s \ = \ \delta m'_1+\epsilon q'_1+\theta p'_1+\rho n'_1 \ > \ 0.$$
En faisant agir le sous-groupe unipotent $H$ de $G$ on se ram\`ene au cas
o\`u \
\m{P'_1 \ = \ X\oplus M'_1},
\m{X} \'etant un sous-espace vectoriel de \ \m{\Hom(\kf_1,\kf_2)\ot M_2}.
En changeant \'eventuellement le rel\`evement de \m{f_2} servant \`a
d\'efinir \m{\phi_1}, on se ram\`ene au cas o\`u \
\m{Q'_1 \ = \ Y\oplus L'_1},
\m{Y} \'etant un sous-espace vectoriel de \ \m{\Hom(\ke_1,\kh_1)\ot M_2}.
Soit \m{M''_2} le plus petit sous-espace vectoriel de \m{M_2} tel que
$$Y\subset \Hom(\kf_1,\kf_2)\ot M''_2.$$

Montrons qu'on peut aussi supposer que \m{M''_2} est le plus petit
sous-espace vectoriel de \m{M_2} tel que
$$X\subset \Hom(\ke_1,\kh_1)\ot M''_2.$$
D'abord on a \ \m{X\subset \Hom(\ke_1,\kh_1)\ot M''_2} : cela d\'ecoule du fait
que pour toute droite \m{L} de \m{\Hom(\ke_1,\kh_1)}, la
restriction de \m{\tau'} \`a \ \m{\Hom(\ke_1,\kf_1)^*\ot L} \ est non nulle
(ceci r\'esultant du fait que \ \m{\kh_1\subset\kf_1\ot\Hom(\kh_1,\kf_1)^*}.
D'autre part on peut, puisque \ \m{\epsilon=\lambda_1>0}, remplacer \m{X}
par \ \m{{\psi'_{11}}^{-1}(Y)\cap(\Hom(\ke_1,\kh_1)\ot M''_2)}, ce qui
assure la minimalit\'e de \m{M''_2}.

Soit \
\m{N''_1 \ = \ \phi((\Hom(\kf_1,\kf_2)\ot M''_2)\oplus M'_1)}.
Alors on a
$$(f_1,f_2)(\ke_1\ot L'_1)\subset (\kf_1\ot M'_1)\oplus(\kf_2\ot M''_2), \ \ \
(g_1,g_2)((\kf_1\ot M'_1)\oplus(\kf_2\ot M''_2))\subset \kg_1\ot N''_1.$$
Posons \ \m{l'_1=\dim(L'_1)}, \m{m'_1=\dim(M'_1)}, \m{n''_1=\dim(N''_1)}, 
\m{m''_2=\dim(M''_2)} \ et
$$t \ = \ \lambda_1l'_1+\mu_1m'_1+\mu_2m''_2+\nu_1n''_1.$$
Alors on a
$$t \ = s + \delta(m''_2-m'_2)
+\epsilon(bm''_2-x)+\theta(am''_2-y)+\rho(n''_1-n'_1).$$
Pour montrer que \m{(*)} n'est pas semi-stable, il suffit donc de montrer que
$$u \ = \ \delta(m''_2-m'_2)+\epsilon(bm''_2-x)+\theta(am''_2-y)
+\rho(n''_1-n'_1) \ \geq \ 0.$$
C'est \'evident si \ \m{am''_2-y=0}, car alors \ \m{n''_1-n'_1=0}. On 
peut donc supposer que \ \m{am''_2-y>0}, ce qui entra\^ine \
\m{c_2(m''_2)>0}. On a alors
$$m''_2-m'_2\geq\frac{am''_2-y}{a}, \ \
bm''_2-x\geq\frac{am''_2-y}{c_2(m''_2)}, \ \ 
n''_1-n'_1\leq c_1(m''_2)(am''_2-y).$$
Donc
$$u \ \geq (am''_2-y)(\frac{\delta}{a}+\frac{\epsilon}{c_2(m''_2)}+\theta
+\rho c_1(m''_2)),$$
c'est-\`a-dire
$$u \ \geq (am''_2-y)(\frac{\mu_2}{a}-(\frac{b}{a}-\frac{1}{c_2(m''_2)})\lambda_1
+\nu_1c_1(m''_2)).$$
On en d\'eduit la proposition. \end{proof}

\bigskip

{\bf Remarque : } En g\'en\'eral, on a
$$b-\frac{a}{c_2(m_2)} \ \leq \ 0.$$
C'est le cas par exemple si \m{\tau'} est stable pour l'action du groupe
r\'eductif \hfil\break
\ \m{SL(\Hom(\ke_1,\kh_1))\times SL(\Hom(\kf_1,\kf_2))}.

\newpage

\begin{sub}{\bf Projectivit\'e des vari\'et\'es de modules}\end{sub}

\begin{subsub}{D\'efinition d'une constante}\end{subsub}

Soient
$$\tau'' : \Hom(\ke_1,\kf_1)^*\ot\Hom(\kf_1,\kg_1)^*\ot\Hom(\kk_1,\ke_1)^*\lra
\Hom(\kf_2,\kg_1)^*$$
l'application lin\'eaire d\'efinie par : 
$$\tau''=\tau\circ(\tau'\ot I_{\Hom(\kf_1,\kg_1)^*}),$$
et pour tout entier positif \m{k},
$$\tau''_k=\tau''\ot I_{\C^k}.$$
Soit \m{G''_k} l'ensemble des sous-espaces vectoriels propres \m{K} de
\ \m{\Hom(\kk_1,\ke_1)^*\ot\C^{k}} \  tels que pour tout sous-espace
vectoriel propre $V$ de \m{\C^{k}}, $K$ ne soit pas contenu dans 
\m{\Hom(\kk_1,\ke_1)^*\ot V}. On pose
$$c(k)= \supp_{K\in{G''_k}}(\frac{\codim(
\tau''_k(\Hom(\ke_1,\kf_1)^*\ot\Hom(\kf_1,\kg_1)^*\ot K))}{\codim(K)}).$$
On montre ais\'ement qu'on a \
\m{c(k) \ \leq \ c_1(k)c_2(k)}.

\vskip 1cm

\begin{subsub}{Projectivit\'e}\end{subsub}

\begin{xtheo}
Il existe un quasi-bon quotient \ \m{Q^{ss}_C//G}, qui est une vari\'et\'e
projective, si les conditions suivantes sont r\'ealis\'ees :

1 - On a 
$$\mu_2-a\nu_1c_1(m_2)>0, \ \ \mu_2-a\mu_1-b\lambda_1>0.$$

2 - Si \ \m{\mu_1 < 0} \ et \ \m{\lambda_1+\nu_1c_(m_2)+
\mu_1c_2(m_2)<0}, alors
$$\mu_2-a\mu_1+b\nu_1c(m_2)+b\mu_1c_2(m_2)\geq 0.$$
\end{xtheo}

\begin{proof}
Le bon quotient \m{{Q''}^{ss}_C//G} existe et est projectif, car \m{G''}
est r\'eductif. D'apr\`es le th\'eor\`eme 3.8 et la remarque qui suit la
proposition 4.3, le quasi-bon quotient \m{Q^{ss}_C//G} existe et est projectif 
si \m{{Q''}^{ss}_C} est contenu dans l'ouvert 
$U$ de \m{Q''_C} contenant les orbites des mutations des points de \m{Q_C}.
Compte tenu des deux \'etapes de mutations, on va caract\'eriser les
complexes
$$\kk_1\ot M_2 \ \hfl{\psi}{} \ \ke_1\ot Q_1
\ \hfl{\psi'}{} \ \kf_1\ot P_1 \ \hfl{\phi}{} \ \kg_1\ot N_1$$
non contenus dans $U$. Soit 
$$f : \Hom(\kk_1,\ke_1)^*\ot M_2\lra Q_1$$
l'application lin\'eaire d\'eduite de \m{\psi}. L'application lin\'eaire
$$g : \Hom(\kk_1,\ke_1)^*\ot \Hom(\ke_1,\kf_1)^*\ot M_2\lra P_1$$
d\'eduite de \m{\psi} et \m{\psi'} se factorise de la fa\c con suivante :
$$\Hom(\kk_1,\ke_1)^*\ot \Hom(\ke_1,\kf_1)^*\ot M_2 \ \hfl{\tau'}{} \
\Hom(\kf_1,\kf_2)\ot M_2 \ \hfl{f'}{} \ P_1,$$
compte tenu des isomorphismes
$$\Hom(\kk_1,\ke_1)^*\simeq \Hom(\ke_1,\kh_1), \ \ \ 
\Hom(\kf_1,\kf_2)^*\simeq \Hom(\kh_1,\kf_1).$$
Alors, \m{(\psi,\psi',\phi)} n'est pas dans $U$ si et seulement si
une des deux propri\'et\'es suivantes est v\'erifi\'ee:
\begin{itemize} 
\item[--] L'application lin\'eaire $f$ n'est pas injective.
\item[--] L'application lin\'eaire $f$ est injective, et $f'$ n'est pas
injective.
\end{itemize}

\bigskip

Supposons que $f$ n'est pas
injective. Soit \m{M'_2} le plus petit sous-espace vectoriel
de \m{M_2} tel que \
\m{\ker(f) \ \subset \ \Hom(\kk_1,\ke_1)^*\ot M'_2}.
Soient
$$d=\dim(\ker(f)), \ \ p = \dim(\ker(\tau')\ot M'_2 \ + \ 
\Hom(\ke_1,\kf_1)^*\ot\ker(f)), \ \ m'_2=\dim(M'_2).$$
Alors on a
$$bm'_2\dim(\Hom(\ke_1,\kf_1))-p \ = \ 
\codim((\tau'\ot I_{M'_2})(\Hom(\ke_1,\kf_1)^*\ot\ker(f))).$$
Il en d\'ecoule que
$$bm'_2\dim(\Hom(\ke_1,\kf_1))-p \ \leq \ c_2(m'_2)(bm'_2-d).$$
Soit 
$$F : \Hom(\ke_1,\kf_1)^*\ot Q_1\lra P_1$$
l'application lin\'eaire d\'eduite de \m{\psi'}, et
$$G : \Hom(\kf_1,\kg_1)^*\ot P_1\lra N_1$$
l'application lin\'eaire d\'eduite de \m{\phi}. Posons
$$Q'_1=\imm(f), \ \ P'_1=F(\Hom(\ke_1,\kf_1)^*\ot Q'_1).$$
Alors on a \
\m{\dim(Q'_1)=bm'_2\dim(\Hom(\ke_1,\kf_1))-p}.
Il r\'esulte de la factorisation pr\'ec\'edente de $g$ qu'on a
$$\dim(P'_1)\leq c_2(m'_2)(bm'_2-d).$$
De m\^eme, l'application lin\'eaire
$$\Hom(\ke_1,\kf_1)^*\ot\Hom(\kf_1,\kg_1)^*\ot\Hom(\kk_1,\ke_1)^*\ot M_2\lra
N_1$$
d\'eduite de \m{\psi}, \m{\psi'}, \m{\phi} se factorise par
\ \m{\tau''\ot I_{M_2}}, et il en d\'ecoule que
$$\dim(N'_1)\leq c(m'_2)(bm'_2-d).$$
Posons \
\m{x = \delta\dim(M'_2)+\epsilon\dim(Q'_1)+\theta\dim(P'_1)+\rho\dim(N'_1)}.
Il faut montrer que \ \m{x > 0}, ce qui montrera que \m{(\psi,\psi',\phi)}
n'est pas semi-stable. On a
$$x \geq (\mu_2-a\mu_1-b\lambda_1)m'_2+(\lambda_1+\nu_1c(m'_2))(bm'_2-d)
+\mu_1\dim(P'_1).$$
Il d\'ecoule ais\'ement des conditions 2- et 3-  du th\'eor\`eme que
\ \m{x > 0}. 

Si $f$ est injective, mais pas $f'$, on d\'eduit de m\^eme de la condition
1- du th\'eor\`eme que \m{(\psi,\psi',\phi)} n'est pas semi-stable.
\end{proof}

\vskip 1.5cm

\begin{sub}{\bf Exemple}\end{sub}

Soient \ \m{n\geq 2} \ un entier et $V$ un espace vectoriel de dimension
n+1. On va \'etudier des complexes sur \ \m{\P_n=\P(V)}, du type
$$\ko(-2)\lra (\ko(-1)\ot M_1)\oplus\ko\lra\ko(1),$$
o\`u $M_1$ est un espace vectoriel tel que \ \m{0<\dim(M_1)\leq n+1}.

On a
$$c_1(1)=0, \ \ c_2(1)=\frac{2}{n+2}, \ \ c(1)=0.$$
Compte tenu du th\'eor\`eme 4.4, en supposant que \ \m{\lambda_1=1}, on
obtient un quasi-bon quotient projectif d\`es que \m{\mu_1}, \m{\mu_2} sont
positifs, \m{\nu_1} n\'egatif, et
$$\mu_2>(n+1)\mu_1+\frac{n(n+1)}{2}.$$
Pour ces valeurs de la polarisation, les quotients sont tous les m\^emes,
On n'obtient donc dans ce cas qu'un seul quasi-bon quotient projectif. 

Traitons plus pr\'ecis\'ement le cas de \m{\P_2}, avec \ \m{m_1=3}.
On montre ais\'ement que \m{Q_C} est irr\'eductible. En cas d'existence
de points stables, les quotients seront donc des vari\'et\'es 
irr\'eductibles de dimension 6.

Un complexe
$$\ko(-2)\lra(\ko(-1)\ot\C^{3})\oplus\ko\lra\ko(1)$$
\'equivaut \`a une paire 
$$((z_1,z_2,z_3,q_0),(q_1,q_2,q_3,z_0)),$$ 
o\`u \ \m{z_i\in V^*}, \m{q_i\in S^2V^*}, et 
$$q_0z_0+q_1z_1+q_2z_2+q_3z_3=0.$$
L'ouvert de \m{Q_C} o\`u \m{z_1}, \m{z_2}, \m{z_3} sont lin\'eairement
ind\'ependants est non vide. Pour un complexe dans cet ouvert, on peut se
ramener, en faisant agir le sous-groupe unipotent maximal, au cas o\`u
\ \m{q_0=0}. Il en d\'ecoule imm\'ediatement que s'il existe des complexes
stables, on doit avoir
$$\mu_2 \ > \ 0.$$
Il n'existe en fait que trois quotients projectifs distincts, qui
correspondent aux cas o\`u \m{\mu_1} est n\'egatif, nul, ou positif.
Dans les deux cas o\`u \m{\mu_1} est non nul, la semi-stabilit\'e d'un 
complexe entra\^ine sa stabilit\'e. 

Le cas o\`u \ \m{\mu_1>0} \ est celui qu'on peut traiter en appliquant
directement le th\'eor\`eme 4.4. Dans ce cas un complexe $x$ est stable 
seulement si le morphisme de gauche est non nul, \m{z_0\not = 0}, et si
pour tout complexe dans la $H$-orbite de $x$, d\'efini par la paire
$$((z_1,z_2,z_3,q'_0),(q'_1,q'_2,q'_3,z_0)),$$ 
\m{q'_1}, \m{q'_2} et \m{q'_3} sont lin\'eairement ind\'ependants.
Le quotient existe et est projectif. Il est non vide, car si
\m{(z_1,z_2,z_3)} est une base de \m{V^*}, le complexe d\'efini par
$$((z_1,z_2,z_3),(z_3^2,-z_2z_3,z_2^2-z_1z_3,z_1))$$
est stable. 

Le cas o\`u \ \m{\mu_1<0} \ ne peut pas \^etre trait\'e directement
(mais on peut le faire en consid\'erant les complexes duaux). Dans ce cas,
un complexe d\'efini par
$$((z_1,z_2,z_3,q_0),(q_1,q_2,q_3,z_0))$$
est stable si et seulement si \m{(z_1,z_2,z_3)} est une base de $V^*$. 
On peut donner une description compl\`ete du quotient, qui est isomorphe
\`a la vari\'et\'e $X$ suivante : soient \m{(z_1,z_2,z_3)} une base de 
\m{V^*}, $E$ le sous-espace vectoriel de \ \m{S^2V^*\ot\C^{3}} \ 
constitu\'e des triplets \m{(q_1,q_2,q_3)} tels que
$$z_1q_1+z_2q_2+z_3q_3=0,$$
$H'$ le sous-espace vectoriel de \ \m{V^*\ot\C^{3}} \ constitu\'e des
triplets \m{(\phi_1,\phi_2,\phi_3)} tels que
$$\phi_1z_1+\phi_2z_2+\phi_3z_3=0.$$
On a un morphisme injectif de fibr\'es vectoriels sur \ 
\m{\P_2=\P(V)} :
$$\Phi : \ko(-1)\ot H'\lra\ko\ot E$$
$$z_0\ot(\phi_1,\phi_2,\phi_3)\longmapsto (z_0\phi_1,z_0\phi_2,z_0\phi_3)$$
On prend alors
$$X \ = \ \P(\coker(\Phi)).$$
 
\vskip 2.5cm

\section{Vari\'et\'es de modules de morphismes}

Soient \m{\ke_1}, \m{\ke_2}, \m{\kf_1} des faisceaux coh\'erents sur une 
vari\'et\'e projective $X$, et \m{M_1}, \m{M_2}, \m{N_1} des espaces
vectoriels de dimension finie. On pose
$$a=\dim(\Hom(\ke_1,\ke_2)), \ \ m_1=\dim(M_1), \ \ m_2=\dim(M_2),
\ n_1=\dim(N_1).$$
On s'int\'eresse \`a des morphismes du
type 
$$(*) \ \ \ \ \ \ \ 
(\ke_1\ot M_1)\oplus(\ke_2\ot M_2)\ \hfl{(f_1,f_2)}{} \ \kf_1\ot N_1.$$
On fait les hypoth\`eses suivantes :
\begin{itemize}
\item[--] Les faisceaux $\ke_1$, $\ke_2$, $\kf_1$ sont simples, et 
$$\Hom(\ke_2,\ke_1)=\Hom(\kf_1,\ke_1)=\Hom(\kf_1,\ke_2)=\nsp.$$
\item[--] Le morphisme canonique
$$\ke_1\ot\Hom(\ke_1,\ke_2)\lra\ke_2$$
est surjectif. On note $\kh_1$ son noyau. 
\item[--] On a \ $\Ext^1(\ke_2,\ke_1)=\nsp$, ce qui entra\^ine un
isomorphisme canonique 
$$\Hom(\kh_1,\ke_1)\simeq\Hom(\ke_1,\ke_2)^*.$$
\end{itemize}

\bigskip

On se trouve dans la situation de la proposition 3.1, avec \ \m{\ke=0},
\m{\Gamma=\ke_1}, \m{M=M_1}, \m{\kg=\ke_2\ot M_2}, \m{\kf=\kf_1\ot N_1}. On a
dans ce cas
$$Q_C=W_C=\Hom((\ke_1\ot M_1)\oplus(\ke_2\ot M_2),\kf_1\ot N_1).$$
Le groupe $G$ op\'erant sur $W_C$ est
$$G \ = \ \Aut((\ke_1\ot M_1)\oplus(\ke_2\ot M_2))\times GL(N_1).$$
On peut voit \ \m{\Aut((\ke_1\ot M_1)\oplus(\ke_2\ot M_2))} \ comme 
constitu\'e de matrices du type
$$\left(\begin{array}{cc}g_1 & 0 \cr \phi & g_2\cr\end{array}\right),$$
avec \ \m{g_1\in GL(M_1)}, \m{g_2\in GL(M_2)} \ et \ 
\m{\phi\in\Hom(\Hom(\ke_1,\ke_2)^*\ot M_1,M_2)}.
Les mutations des morphismes \m{(*)} sont d\'ecrites au \para 5.2. On les
appelle des {\em mutations directes} car elles nous 
ram\`enent imm\'ediatement
\`a une action d'un groupe r\'eductif, et permettent de construire les
vari\'et\'es de modules de morphismes. C'est la m\'ethode appliqu\'ee dans
\cite{dr_tr} (les r\'esultats sont rappel\'es dans le \para 5.2).

On peut appliquer la proposition 3.1 d'une autre mani\`ere. C'est 
simple \`a voir si \m{\ke_1}, \m{\ke_2} et \m{\kf_1} sont localement libres :
les complexes de type \m{(*)} sont \'equivalents \`a des complexes du
type
$$\kf_1^*\ot N_1^*\lra (\ke_2^*\ot M_2^*)\oplus(\ke_1^*\ot M_1^*),$$
et on applique la proposition 2.1 en prenant \ \m{\ke=\kf_1^*\ot N_1^*},
\m{\Gamma=\ke_2}, \m{M=M_2^*}, 
\m{\kg=\ke_1^*\ot M_1^*}, \m{\kf=0}. Les
mutations obtenues sont des morphismes de type (*) (avec d'autres 
faisceaux et d'autres espaces vectoriels). On peut ensuite appliquer la
mutation directe \`a ces morphismes, et on obtient de nouveaux quotients
(cf. \para 5.3).

\vskip 1.5cm

\begin{sub}{\bf Notions de (semi-)stabilit\'e et vari\'et\'es de modules de
morphismes}\end{sub}

Le groupe $G$ poss\`ede deux sous-groupes importants : le premier est le
sous-groupe normal unipotent maximal \'evident $H$, isomorphe au groupe
additif 
\m{\Hom(\Hom(\ke_1,\ke_2)^*\ot M_1,M_2))}. Le second est le sous-groupe
r\'eductif
$$G_{red} \ = \ GL(M_1)\times GL(M_2)\times GL(N_1),$$
dont l'inclusion dans $G$ induit un isomorphisme \ 
\m{G_{red}\simeq G/H}.

L'action de \m{G_{red}} est un cas particulier des actions \'etudi\'ees 
dans \cite{king}. Soient \m{\lambda_1}, \m{\lambda_2}, \m{\mu_1}
des nombres rationnels non nuls tels que
$$\lambda_1\dim(M_1)+\lambda_2\dim(M_2)-\mu_1\dim(N_1)=0.$$ 

\bigskip

\begin{defin}
Un point \m{(f_1,f_2)} de $W_C$ est dit {\em $G_{red}$-semi-stable} (resp.
{\em $G_{red}$-stable}) relativement \`a \m{(\lambda_1,\lambda_2,\mu_1)} si
pour tous sous-espaces vectoriels \ \m{M'_1\subset M_1}, 
\m{M'_2\subset M_2}, \m{N'_1\subset N_1}, avec \
\m{(M'_1,M'_2,N'_1)\not = (M_1,M_2,N_1)} \ ou \m{(\nsp,\nsp,\nsp)}, tels que
$$f_1(\ke_1\ot M'_1)\subset \kf_1\ot N'_1, \ \ \
f_2(\ke_2\ot M'_2)\subset \kf_1\ot N'_1,$$
on a
$$\lambda_1\dim(L'_1)+\lambda_2\dim(M'_2)-\mu_1\dim(N'_1)\leq 0
\ \ \ {\rm (resp. \ \ } \ < \ {\rm )}.$$
\end{defin}

\bigskip 

On dit que \m{(\lambda_1, \lambda_2, \mu_1)} est une {\em polarisation}
de l'action de $G$ sur $W_C$. On note \m{W^{ss}_{C,red}} (resp.
\m{W^s_{C,red}}) l'ouvert de \m{W_C} constitu\'e
des points \m{G_{red}}-semi-stables (resp. \m{G_{red}}-stables). D'apr\`es
\cite{king}, il existe un bon quotient \ \m{W^{ss}_{red,C}//G} \ et un
quotient g\'eom\'etrique lisse \m{W^s_{C,red}/G}. Mais ce ne sont pas les
quotients que nous recherchons.

\bigskip

\begin{defin}
Un point \m{x} de $W_C$ est dit {\em $G$-semi-stable} (resp.
{\em $G$-stable}) relativement \`a \m{(\lambda_1,\lambda_2,\mu_1)} si
tous les points de l'orbite \m{H.x} sont $G_{red}$-semi-stables (resp.
$G_{red}$-stables).
\end{defin}

\bigskip 

On note \m{W_C^{ss}} (resp. \m{W_C^s}) l'ouvert de \m{W_C} constitu\'e
des points \m{G}-semi-stables (resp. \m{G}-stables). On cherche
\`a prouver l'existence de bons quotients \ \m{W_C^{ss}//G}. De tels 
quotients seront appel\'es des {\em vari\'et\'es de modules de complexes}.

\medskip

On montre ais\'ement que si \m{W^{s}} est non vide on doit avoir
\ \m{\lambda_1>0, \ \ \lambda_2>0, \ \ \mu_1>0}.
On supposera par la suite que ces in\'egalit\'es sont v\'erifi\'ees.
On peut alors {\em normaliser} la polarisation, c'est-\`a-dire supposer
que
$$\lambda_1\dim(M_1)+\lambda_2\dim(M_2)=1, \ \ \mu_1=\frac{1}{\dim(M_1)}.$$

\vskip 1.5cm

\begin{sub}{\bf Construction des vari\'et\'es de modules par les
mutations directes}\end{sub}

La mutation de \m{(*)} est un complexe
$$(**) \ \ \ \ \ \ \ \
\kh_1\ot M_2\ \hfl{\phi_1}{} \ \ke_1\ot P_1\ \hfl{\phi_2}{} \ \kf_1\ot N_1,$$
avec
$$P_1 \ = \ (\Hom(\ke_1,\ke_2)\ot M_2)\oplus M_1.$$
Le morphisme \m{\phi_1} est d\'efini par l'inclusion
$$\Hom(\kh_1,\ke_1)^*\ot M_2=\Hom(\ke_1,\ke_2)\ot M_2 \ \subset \ P_1.$$
Le morphisme \m{\phi_2} provient d'une application lin\'eaire
$$F : \Hom(\ke_1,\kf_1)^*\ot P_1\lra N_1,$$
qui est \'egale \`a celle d\'eduite de \m{f_1} sur \
\m{\Hom(\ke_1,\kf_1)^*\ot M_1}. Sur l'autre facteur, $F$ est la compos\'ee
$$\Hom(\ke_1,\kf_1)^*\ot\Hom(\ke_1,\ke_2)\ot M_2\lra\Hom(\ke_2,\kf_1)^*\ot M_2
\lra N_1,$$
la premi\`ere application provenant de la composition
$$\sigma : \Hom(\ke_1,\ke_2)\ot\Hom(\ke_2,\kf_1)\lra\Hom(\ke_1,\kf_1),$$
et la seconde de \m{f_2}.

Le groupe op\'erant sur la vari\'et\'e \m{Q'_C} des complexes \m{(**)} est
$$G' \ = \ GL(M_2)\times GL(P_1)\times GL(N_1),$$
qui est r\'eductif. Une notion de (semi-)stabilit\'e pour les points de
\m{Q'_C} est d\'efinie par une suite \m{(\alpha,\beta,\gamma)} de nombres
rationnels non nuls telle que
$$\alpha\dim(M_2)+\beta\dim(P_1)+\gamma\dim(N_1) = 0.$$
Un complexe \m{(\phi'_1,\phi'_2)} de type \m{(**)}
est semi-stable (resp.stable) relativement \`a \m{(\alpha,\beta,\gamma)} si
et seulement si pour tous sous-espaces vectoriels
$$M'_2\subset M_2, \ \ P'_1\subset P_1, \ \ N'_1\subset N_1,$$
avec \ \m{(M'_2,P'_1,N'_1)\not = (M_2,P_1,N_1)} \ ou \m{(\nsp,\nsp,\nsp)}, 
tels que
$$\phi'_1(\kh_1\ot M'_2)\subset\ke_1\ot P'_1, \ \
\phi'_2(\ke_1\ot P'_1)\subset\kf_1\ot N'_1,$$
on a
$$\alpha\dim(M'_2)+\beta\dim(P'_1)+\gamma\dim(N'_1)\leq 0 \ \ \
\ {\rm (resp. \ } < \ {\rm )}.$$ 
S'il existe des complexes stables on doit avoir
\ \m{\alpha > 0, \ \ \gamma>0}.

Soient \ \m{M'_1\subset M_1}, \m{M'_2\subset M_2}, \m{N'_1\subset N_1} \
des sous-espaces vectoriels tels que
$$f_1(\ke_1\ot M'_1)\subset\kf_1\ot N'_1, \ \
f_2(\ke_2\ot M'_2)\subset\kf_1\ot N'_1.$$
On pose
$$P'_1 \ = \ (\Hom(\ke_1,\ke_2)\ot M'_2)\oplus M'_1.$$
On a alors
$$\phi_1(\kh_1\ot M'_2)\subset\ke_1\ot P'_1, \ \ 
\phi_2(\ke_1\ot P'_1)\subset\kf_1\ot N'_1.$$

On suppose que
$$\alpha=\lambda_2-a\lambda_1, \ \beta=\lambda_1, \ \gamma=-\mu_1.$$
On a alors
$$\alpha\dim(M'_2)+\beta\dim(P'_1)+\gamma\dim(N'_1) \ = \
\lambda_1\dim(M'_1)+\lambda_2\dim(M'_2)-\mu_1\dim(N'_1).$$
On en d\'eduit imm\'ediatement la

\bigskip

\begin{xprop}
Si le complexe \m{(**)} est semi-stable (resp. stable) relativement \`a 
\m{(\alpha,\beta,\gamma)}, le complexe (*) est $G$-semi-stable (resp.
$G$-stable) relativement \`a \m{(\lambda_1,\lambda_2,\mu_1)}.
\end{xprop}

\bigskip

Pour continuer on utilise une constante analogue \`a celles qui ont \'et\'e
d\'efinies au \para 4.3.1. On consid\`ere l'application canonique
$$\tau : \Hom(\ke_1,\kf_1)^*\ot\Hom(\ke_1,\ke_2)\lra\Hom(\ke_2,\kf_1)^*.$$
Pour tout entier positif $k$, soit
$$\tau_k = \tau\ot I_{\C^k} :
\Hom(\ke_1,\kf_1)^*\ot(\Hom(\ke_1,\ke_2)\ot\C^{k})
\lra\Hom(\ke_2,\kf_1)^*\ot\C^{k}.$$
Soit \m{G_k} l'ensemble des sous-espaces vectoriels propres \m{K} de 
\m{\Hom(\ke_1,\ke_2)\ot\C^{k}} tels que pour tout sous-espace vectoriel
propre \ \m{V\subset\C^{k}}, \m{K} ne soit pas contenu dans

\noindent\m{\Hom(\ke_1,\ke_2)\ot V}. On pose
$$c_0(k)= \supp_{K\in{G_k}}(\frac{\codim(
\tau_k(\Hom(\ke_1,\kf_1)^*\ot K))}{\codim(K)}).$$
On d\'emontre dans \cite{dr_tr} le

\bigskip

\begin{xtheo}
Si
$$\frac{\lambda_2}{\lambda_1}>a, \ \ \ \lambda_2>\frac{a}{n_1}c_0(m_2),$$ 
il existe un bon quotient projectif \ \m{W_C^{ss}//G}, et un quotient
g\'eom\'etrique \ \m{W_C^{s}/G}, qui est un ouvert de \ \m{W_C^{ss}//G}.
\end{xtheo}

\bigskip

On montre en fait que sous les hypoth\`eses du th\'eor\`eme, la 
r\'eciproque de la proposition 5.1 est vraie. On peut conclure en 
utilisant le th\'eor\`eme 3.8.

\newpage

\begin{sub}{\bf Mutations indirectes}\end{sub}

\begin{subsub}{D\'efinition et premi\`eres propri\'et\'es}\end{subsub}

On utilise maintenant la deuxi\`eme mani\`ere d'effectuer des mutations
constructives.
Comme indiqu\'e au d\'ebut du \para 5, c'est plus facile \`a
voir si \m{\ke_1}, \m{\ke_2} et \m{\kf_1} sont localement libres. Mais les
constructions \'etant purement formelles, il n'est pas n\'ecessaire de faire
cette supposition. On fait les hypoth\`eses suppl\'ementaires suivantes :
\begin{itemize}
\item[--] Le morphisme canonique
$$\ke_1\lra\ke_2\ot\Hom(\ke_1,\ke_2)^*$$
est injectif. Soit $\kg_1$ son conoyau.
\item[---] La composition
$$\sigma : \Hom(\ke_1,\ke_2)\ot\Hom(\ke_2,\kf_1)\lra\Hom(\ke_1,\kf_1)$$
est surjective.
\item[--] On a \ $\Hom(\kf_1,\kg_1)=\nsp$.
\end{itemize}

\bigskip

Cette derni\`ere hypoth\`ese n'est pas strictement indispensable,
on le verra plus loin.
On a, puisque \ \m{\Ext^1(\ke_2,\ke_1)=\nsp}, un
isomorphisme canonique \
\m{\Hom(\ke_2,\kg_1) \ \simeq \ \Hom(\ke_1,\ke_2)^*}.

Une mutation d'un morphisme \m{(*)} est un morphisme
$$(***) \ \ \ \ \ \ \ 
\ke_2\ot Q_1 \ \hfl{(\psi_1,\psi_2)}{} \ (\kg_1\ot M_1)\oplus(\kf_1\ot N_1),$$
avec
$$Q_1 \ = \  (\Hom(\ke_1,\ke_2)^*\ot M_1)\oplus M_2.$$
Le morphisme \m{\psi_1} est nul sur \ \m{\ke_2\ot M_2}, et sur l'autre facteur
c'est 
$$ev\ot I_{M_1} : \ke_2\ot\Hom(\ke_1,\ke_2)^*\ot M_1 = 
\ke_2\ot\Hom(\ke_2,\kg_1)\ot M_1\lra\kg_1\ot M_1,$$
\m{ev} d\'esignant le morphisme d'\'evaluation. Le morphisme \m{\psi_2} est
\'egal \`a \m{f_2} sur 

\noindent \m{\Hom(\ke_2,\kf_1)^*\ot M_2}, et sur l'autre facteur
c'est un rel\`evement de \m{f_1}, qui existe \`a cause de la seconde
hypoth\`ese suppl\'ementaire. Remarquons que \m{\psi_2} est d\'efini \`a un
\'el\'ement pr\`es de \m{\Hom(\kg_1\ot M_1, \kf_1\ot N_1)}.

Soit \m{W'_C} l'espace vectoriel des morphismes de type \m{(***)}, \m{G'} le
groupe op\'erant sur \m{W'_C} :
$$G' \ = \ GL(Q_1)\times\Aut((\kg_1\ot M_1)\oplus(\kf_1\ot N_1)).$$
Si on omet la derni\`ere hypoth\`ese suppl\'ementaire, il faut remplacer
\m{G'} par un groupe plus petit. Soient \m{H'} le sous-groupe normal unipotent
maximal de \m{G'}, isomorphe au groupe additif
\m{\Hom(\kg_1\ot M_1, \kf_1\ot N_1)}, et \m{G'_{red}} le sous-groupe r\'eductif :
$$G'_{red} \ = \  GL(Q_1)\times GL(M_1)\times GL(N_1).$$

Soient \m{\nu_1}, \m{\nu_2} des nombres rationnels positifs tels que
$$\nu_1\dim(N_1)+\nu_2\dim(M_1) = 1.$$
Alors \m{(1/\dim(Q_1),\nu_2,\nu_1)} d\'efinit une notion de semi-stabilit\'e
pour l'action de \m{G'} sur \m{W'_C}. Rappelons qu'un point 
\m{(\psi'_1,\psi'_2)}
de \m{W'_C} est {\em \m{G'_{red}}-semi-stable} 
(resp. {\em \m{G'_{red}}-stable})
relativement \`a \m{(1/\dim(Q_1),\nu_2,\nu_1)} si et seulement si pour tous 
sous-espaces vectoriels 
$$Q'_1\subset Q_1, \ \ M'_1\subset M_1, \ \ N'_1\subset N_1,$$
avec \ \m{(Q'_1,M'_1,N'_1)\not =  (Q_1,M_1,N_1)} \ ou \m{(\nsp,\nsp,\nsp)},
tels que
$$\psi'_1(\ke_2\ot Q'_1)\subset\kg_1\ot M'_1, \ \ 
\psi'_2(\ke_2\ot Q'_1)\subset\kf_1\ot N'_1,$$
on a
$$\frac{\dim(Q'_1)}{\dim(Q_1)} - \nu_2\dim(M'_1)-\nu_1\dim(N'_1)\leq 0,
\ \ \ \ {\rm (resp. \ } \ < \ {\rm)}.$$
On dit que \m{(\psi'_1,\psi'_2)} est {\em \m{G'}-semi-stable} (resp.
{\em \m{G'}-stable}) si tous les points de sa \m{H'}-orbite sont
\m{G'_{red}}-semi-stables (resp. \m{G'_{red}}-stables).

Soient \ \m{M'_1\subset M_1}, \m{M'_2\subset M_2}, \m{N'_1\subset N_1} \ des
sous-espaces vectoriels tels que
$$f_1(\ke_1\ot M'_1)\subset\kf_1\ot N'_1, \ \ 
f_2(\ke_2\ot M'_2)\subset\kf_1\ot N'_1.$$
On pose 
$$Q'_1 \ = \ (\Hom(\ke_1,\ke_2)^*\ot M'_1)\oplus M'_2.$$
Il est clair qu'on peut choisir la mutation \m{(***)} de telle sorte que
$$\psi_1(\ke_2\ot Q'_1)\subset\kg_1\ot M'_1, \ \
\psi_2(\ke_2\ot Q'_1)\subset\kf_1\ot N'_1.$$
Posons
$$\nu_1=\frac{1}{\dim(Q_1)n_1\lambda_2}, \ \
\nu_2=\frac{a\lambda_2-\lambda_1}{\dim(Q_1)\lambda_2},$$
de telle sorte qu'on a
$$\frac{\dim(Q'_1)}{\dim(Q_1)} - \nu_2\dim(M'_1)-\nu_1\dim(N'_1)=
\frac{\lambda_1\dim(M'_1)+\lambda_2\dim(M'_2)-\dim(N'_1)/n_1}
{\lambda_2\dim(Q_1)}.$$
Alors on a d'apr\`es ce qui pr\'ec\`ede la

\bigskip

\begin{xprop}
Si le morphisme \m{(***)} est \m{G'}-semi-stable (resp.  \m{G'}-stable)
relativement \`a \m{(1/\dim(Q_1),\nu_2,\nu_1)}, alors le morphisme \m{(*)} est
\m{G}-semi-stable (resp.  \m{G}-stable) relativement \`a 
\m{(\lambda_1,\lambda_2,1/\dim(N_1))}.
\end{xprop}

\bigskip

S'il existe des morphismes \m{G'}-stables dans \m{W'_C}, on a \ \m{\nu_1>0},
c'est-\`a-dire \
\m{a\lambda_2-\lambda_1>0}.
On supposera que c'est le cas dans toute la suite.

\vskip 1cm

\begin{subsub}{Cas d'\'equivalence des (semi-)stabilit\'es}\end{subsub}

Soit
$$\sigma : \Hom(\ke_2,\kf_1)^*\ot\Hom(\ke_1,\ke_2)^*\lra\Hom(\ke_1,\kf_1)^*$$
une application lin\'eaire dont la composition avec la transpos\'ee de la
composition
$$\Hom(\ke_1,\ke_2)\ot\Hom(\ke_2,\kf_1)\lra\Hom(\ke_1,\kf_1)$$
est l'identit\'e de \m{\Hom(\ke_1,\kf_1)^*}. 
On d\'efinit ici des constantes analogues \`a celles du \para 4.3.1. 
Pour tout entier positif $k$, soit
$$\sigma_k = \sigma\ot I_{\C^k} :
\Hom(\ke_2,\kf_1)^*\ot(\Hom(\ke_1,\ke_2)^*\ot\C^{k})
\lra\Hom(\ke_1,\kf_1)^*\ot\C^{k}.$$
Soit \m{G_k} l'ensemble des sous-espaces vectoriels propres \m{K} de 
\m{\Hom(\ke_1,\ke_2)^*\ot\C^{k}} tels que pour tout sous-espace vectoriel
propre \ \m{V\subset\C^{k}}, \m{K} ne soit pas contenu dans

\noindent\m{\Hom(\ke_1,\ke_2)^*\ot V}. On pose
$$c'_0(k)= \supp_{K\in{G_k}}(\frac{\codim(
\sigma_k(\Hom(\ke_1,\ke_2)^*\ot K))}{\codim(K)}).$$
Il est clair qu'on a \
\m{c'_0(k+1)\geq c'_0(k)}.

\medskip

{\bf Remarque : } les constantes pr\'ec\'edentes d\'ependent aussi
a priori du choix de \m{\sigma}.

\bigskip

\begin{xprop}
On suppose que 
$$\lambda_2\geq \frac{c'_0(m_1)}{n_1}.$$
Si le morphisme \m{(*)} est \m{G}-semi-stable 
relativement \`a \m{(\lambda_1,\lambda_2,1/\dim(N_1))}, alors le morphisme
\m{(***)} est \m{G'}-semi-stable relativement \`a \m{(1/\dim(Q_1),\nu_2,\nu_1)}.
\end{xprop}

\bigskip

\begin{proof}
Soient \ \m{Q'_1\subset Q_1}, \m{M'_1\subset M_1}, \m{N'_1\subset N_1} \ des
sous-espaces vectoriels tels que
$$\psi_1(\ke_2\ot Q'_1)\subset \kg_1\ot M'_1, \ \ 
\psi_2(\ke_2\ot Q'_1)\subset \kf_1\ot N'_1$$
et
$$X = \frac{\dim(Q'_1)}{\dim(Q_1)}-\nu_2\dim(M'_1)-\nu_1\dim(N'_1)>0.$$
En faisant agir le sous-groupe unipotent $H$ de $G$, on peut supposer que
\m{Q'_1} est de la forme \
\m{Q'_1 \ = \ K\oplus M'_2},
avec \ \m{K\subset\Hom(\ke_1,\ke_2)^*\ot M_1} \ et \ \m{M'_2\subset M_2}.
On peut supposer que \m{M'_1} est le plus petit sous-espace vectoriel de \m{M_1}
tel que \ \m{K\subset\Hom(\ke_1,\ke_2)^*\ot M'_1}. Soit \m{N''_1} l'image de
\ \m{\Hom(\ke_1,\kf_1)^*\ot M'_1} \ par l'application lin\'eaire
$$\Hom(\ke_1,\kf_1)^*\ot M_1\lra N_1$$
d\'eduite de \m{f_1}. Alors on a
$$\dim(N''_1)-\dim(N'_1)\leq c'_0(m_1)(a\dim(M'_1)-\dim(K)).$$
On a
$$f_1(\ke_1\ot M'_1)\subset \kf_1\ot N''_1, \ \ \ 
f_2(\ke_2\ot M'_2)\subset \kf_1\ot N''_1.$$
On va montrer que
$$X'= \lambda_1\dim(M'_1)+\lambda_2\dim(M'_2)-\frac{\dim(N''_1)}{n_1} > 0,$$
ce qui prouvera que \m{(*)} n'est pas semi-stable. On a
$$\frac{X'}{\dim(Q_1)\lambda_2} \ = \ X + \frac{a\dim(M'_1)-\dim(K)}{\dim(Q_1)}
-\nu_1(\dim(N''_1)-\dim(N'_1)),$$
donc
$$\frac{X'}{\dim(Q_1)\lambda_2} \ > \ 
(\frac{1}{\dim(Q_1)}-\nu_1c'_0(m_1))(a\dim(M'_1)-\dim(K)),$$
c'est-\`a-dire
$$X' \ > \ (\lambda_2-\frac{c'_0(m_1)}{n_1})(a\dim(M'_1)-\dim(K)) \geq 0.$$
\end{proof}

\vskip 1cm

\begin{subsub}{Egalit\'e des quotients et projectivit\'e}\end{subsub}

On se place dans l'hypoth\`ese du th\'eor\`eme 5.5. On peut se poser la question
de l'\'egalit\'e des quotients \ \m{{W'_C}^{ss}//G'} \ et \ \m{W_C^{ss}//G}.
Soit \m{W'_0} l'ouvert de \m{W'_C} constitu\'e des morphismes 
\m{(\psi'_1,\psi'_2)} tels que \m{\psi'_1} induise une surjection
$$F : Q_1\lra\Hom(\ke_1,\ke_2)^*\ot M_1.$$
Cet ouvert est donc constitu\'e des morphismes dont un \'el\'ement de la 
$G'$-orbite peut \^etre obtenu comme mutation d'un \'el\'ement de $W_C$.
Si \ \m{(\psi'_1,\psi'_2)\in W'_C\backslash W'_0}, il existe un sous-espace
vectoriel \ \m{Q'_1\subset Q_1} \ de dimension \ \m{m_2+1} \ contenu dans
\m{\ker(F)}. On a
$$\psi'_1(\ke_2\ot Q'_1)=0.$$
Il en d\'ecoule que si \m{(\psi'_1,\psi'_2)} est \m{G'}-semi-stable 
relativement \`a \m{(1/\dim(Q_1,),\nu_2,\nu_1)}, on doit avoir
$$\frac{m_2+1}{\dim(Q_1)}-\nu_1n_1\geq 0,$$
c'est-\`a-dire
$$\lambda_2 \ \leq \frac{1}{m_2+1}.$$
On en d\'eduit le

\bigskip

\begin{xtheo}
On suppose que
$$\lambda_2\geq \frac{c'_0(m_1)}{n_1} \ \ \ \ {\rm et} \ \ \ \
\lambda_2>\frac{1}{m_2+1}.$$
S'il existe un quasi-bon quotient \ \m{{W'_C}^{ss}//G'} \ pour la polarisation
\m{(1/\dim(Q_1,),\nu_2,\nu_1)}, il existe un quasi-bon quotient \ 
\m{W_C^{ss}//G} \ pour la polarisation
\m{(\lambda_1,\lambda_2,1/n_1)}, et il est isomorphe \`a \ \m{{W'_C}^{ss}//G'}.
\end{xtheo}

\vskip 1cm

\begin{subsub}{Cas d'existence d'un quasi-bon quotient projectif}\end{subsub}

On fait ici la synth\`ese du th\'eor\`eme 5.2, du th\'eor\`eme 7.6 de 
\cite{dr2} et des r\'esultats pr\'ec\'edents. 

Pour pouvoir employer le th\'eor\`eme 7.6 de \cite{dr2}, il faut faire les
hypoth\`eses suppl\'ementaires suivantes :
\begin{itemize}
\item[--] Les morphismes
$$\coeva{\ke_1}{\kf_1}, \ \ \coeva{\ke_2}{\kf_1}$$
sont injectifs. On note \ $\kh_1$, $\kh_2$ leurs conoyaux respectifs.
\item[--] On a \ \ $\Ext^1(\kf_1,\ke_1)=\Ext^1(\kf_1,\ke_2)=\nsp$.
\end{itemize}

\bigskip

Soit
$$(\ke_1\ot M_1)\oplus(\ke_2\ot M_2)\lra\kf_1\ot N_1$$
un morphisme tel que l'application lin\'eaire associ\'ee 
$$\phi : (\Hom(\ke_1,\kf_1)^*\ot M_1)\oplus(\Hom(\ke_2,\kf_1)^*\ot M_2)\lra N_1$$
soit surjective. On d\'efinit dans \cite{dr2} un autre type de mutation 
associant au morphisme pr\'ec\'edent un morphisme
$$\kf_1\ot\ker(\phi)\lra(\kh_1\ot M_1)\oplus(\kh_2\ot M_2).$$
En appliquant le th\'eor\`eme 5.2 \`a ces nouveaux morphismes, on obtient de
nouveaux cas d'existence de quasi-bons quotients. 

On aura besoin de deux types de constantes suppl\'ementaires. Soit $k$ un
entier positif. On note \m{c_1(k)} la constante analogue \`a \m{c_0(k)},
obtenue en consid\'erant l'application lin\'eaire
$$\Hom(\ke_2,\kf_1)\ot(\Hom(\ke_1,\ke_2)\ot\C^{k})\lra \Hom(\ke_1,\kf_1)\ot\C^{k}$$
au lieu de \m{\tau_k} (cf. \para 5.2). On note \m{c_2(k)} la constante analogue
 \`a \m{c_0(k)}, obtenue en consid\'erant l'application lin\'eaire
$$\Hom(\ke_2,\kf_1)^*\ot(\ker(\sigma)\ot\C^{k})\lra\Hom(\ke_1,\ke_2)\ot\C^{k}.$$
On pose
$$h_{11}=\dim(\Hom(\ke_1,\kf_1)), \ \ h_{12}=\dim(\Hom(\ke_2,\kf_1)), \ \ 
a'=ah_{12}-h_{11}.$$

Toutes m\'ethodes confondues, on obtient le

\bigskip

\begin{xtheo}
Il existe un quasi-bon quotient projectif \ \m{W^{ss}_C//G} \ dans les cas
suivants : 

1 - On a 
$$\frac{\lambda_2}{\lambda_1} > a, \ \ \ \ \lambda_2 > \frac{a}{n_1}c_0(m_2).$$

2 - On a
$$\lambda_1 < \frac{h_{11}}{n_1}, \ \ \ \ \lambda_2 < \frac{h_{12}}{n_1}, \ \ \ \ 
a\lambda_2-\lambda_1 > \frac{a'}{n_1}, \ \ h_{11}-\lambda_1n_1\geq c_1(m_1)a.$$

3 - On a
$$\lambda_2\geq\frac{c'_0(m_1)}{n_1}, \ \ \ \ \lambda_2 > \frac{1}{m_2+1}, \ \ \ \
a\lambda_2-\lambda_1 \ > \ a'.{\rm Max}(c_2(m_1)\lambda_2, \frac{1}{n_1}).$$
\end{xtheo}

On peut montrer que dans la situation du th\'eor\`eme 5.6, l'ouvert du
quotient correspondant aux points stables est lisse.

\vskip 1.5cm

\begin{sub}{\bf Exemples}\end{sub}

Si  \m{(\lambda_1,\lambda_2,1/n_1)} est une polarisation de l'action de $G$
sur \m{W_C}, on notera
$$\rho \ = \ \frac{\lambda_2}{\lambda_1}.$$
Notons que la polarisation est enti\`erement d\'etermin\'ee par \m{\rho}, ou
par un des nombres \m{\lambda_1}, \m{\lambda_2}.

Dans les exemples qui vont suivre, on utilise implicitement des calculs de
constantes \m{c_0(k)} ou \m{c'_0(k)} qui proviennent de \cite{dr_tr}, 

\newpage

\begin{subsub}{Exemple 1}\end{subsub}

Soit $n$ un entier, avec \ \m{n\geq 2}. On consid\`ere des morphismes
$$(\ko(-2)\ot\C^{2})\oplus\ko(-1)\lra\ko\ot\C^{n+1}$$
sur \m{\P_n}. Cet exemple a d\'ej\`a \'et\'e trait\'e dans \cite{dr_tr},
pour \ \m{n=2}. L'application directe du th\'eor\`eme 5.2 montre qu'il existe
un bon quotient projectif \ \m{W_C^{ss}//G} \ si 
$$\rho \ > \ n+1.$$ 
En utilisant des mutations indirectes on se ram\`ene \`a des morphismes
$$\ko(-1)\ot\C^{2n+3}\lra(Q(-1)\ot\C^{2})\oplus(\ko\ot\C^{n+1}).$$
L'application des th\'eor\`emes 5.2 et 5.5 permet de montrer qu'il existe
un quasi-bon quotient projectif \ \m{W_C^{ss}//G} \ d\`es que 
$$\rho \ > \ 2 + \frac{2}{n}.$$
Si \ \m{n\geq 4}, on obtient ainsi des vari\'et\'es de modules de morphismes
suppl\'ementaires.

\vskip 1cm

\begin{subsub}{Exemple 2}\end{subsub}

On consid\`ere les morphismes
$$(f_1,f_2) : \ko(-2)\oplus\ko(-1)\lra\ko\ot\C^{n+2}$$
sur $\P_n$. C'est un exemple d\'ej\`a donn\'e dans \cite{dr2}.
On sait construire des bons quotients (en utilisant le th\'eor\`eme 5.2)) 
d\`es que
$$\rho > n+1.$$
Mais dans ce cas le quotient est vide. En effet, il existe toujours un
sous-espace vectoriel \ \m{H\subset\C^{n+2}} \ de dimension \m{n+1} tel 
que \ \m{\imm(f_2)\subset \ko\ot H}. On doit donc avoir, si \m{(f_1,f_2)}
est \m{G}-semi-stable relativement \`a \m{(\lambda_1,\lambda_2,1/(n+2))},
$$\lambda_2 - \frac{n+1}{n+2} \leq 0,$$
c'est-\`a-dire \ \m{\rho\leq n+1}. L'application des th\'eor\`emes 5.5 et
5.2 permet de construire un quasi-bon quotient projectif \ 
\m{W^{ss}_C//G} \ d\`es que
$$\rho \ > 1.$$
On am\'eliore l\'eg\`erement le r\'esultat de \cite{dr2} dans le cas o\`u $n$
est impair. 
Les valeurs {\em singuli\`eres} de \m{\rho} sont par d\'efinition celles
pour lesquelles la \m{G}-semi-stabilit\'e n'implique pas la
\m{G}-stabilit\'e. Ces valeurs sont exactement les nombres
$$\rho_k = \frac{k}{n+2-k}$$
pour \m{1\leq k\leq n+1} . Dans ce cas un morphisme \m{(\phi_1,\phi_2)}
$G$-semi-stable non $G$-stable est construit de la fa\c con suivante :
on consid\`ere un sous-espace vectoriel \ \m{H\subset\C^{n+2}} \ de
dimension \m{k}, et on prend pour \m{\phi_2} un morphisme tel que 
\ \m{\imm(\phi_2)\subset\ko\ot H} \ et que $H$ soit le plus petit 
sous-espace vectoriel ayant cette propri\'et\'e. On prend pour \m{\phi_1}
un morphisme tel que l'application lin\'eaire induite
$$H^0(\ko(2))^*\lra\C^{n+2}$$
soit surjective. Toutes le polarisations telles que \m{\rho} soit situ\'e entre
\m{\rho_k} et \m{\rho_{k+1}} donnent la m\^eme notion de (semi-)stabilit\'e.

Notons \m{M_k} le quotient obtenu pour \ \m{\rho=\rho_k}, et \m{M^0_k} celui
obtenu pour \ \m{\rho_{k-1}<\rho<\rho_k}. On sait donc construire \m{M_k} et
\m{M^0_k} pour \ \m{\lbrack\frac{n}{2}\rbrack+2\leq k\leq n+1}. 
Les quotients \m{M^0_k} construits dans \cite{dr2}, c'est-\`a-dire pour \
\m{\lbrack\frac{n+3}{2}\rbrack+1\leq k\leq n+1}, sont des quotients 
g\'eom\'etriques. 

\vskip 1cm

\begin{subsub}{Exemple 3}\end{subsub}
 
C'est une g\'en\'eralisation de l'exemple pr\'ec\'edent. Soit \ \m{n\geq 2} \
un entier. On consid\`ere sur $\P_n$ les morphismes
$$(f_1,f_2) : \ko(-2)\oplus\ko(-1)\lra\ko\ot N_1.$$

\vskip 0.8cm

\subsubsubsection{5.4.3.1}{Conditions d'existence de points stables}

\medskip

Pour qu'il existe des points stables (pour au moins une polarisation), on doit
avoir
$$n_1=\dim(N_1) \ \leq \ \frac{(n+1)(n+2)}{2}+n.$$

Si \ \m{n_1\geq n+1}, il existe toujours un sous-espace vectoriel \m{N'_1}
de \m{N_1} de dimension \m{n+1} tel que l'image de \m{f_2} soit contenue dans
\ \m{\ko\ot N'_1}. Il en d\'ecoule que s'il existe des points stables
relativement \`a la polarisation d\'efinie par \m{\lambda_2}, on doit avoir
$$\lambda_2 \ < \ \frac{n+1}{n_1}.$$
Si c'est le cas on montre ais\'ement qu'il existe toujours des morphismes 
stables.

Les valeurs de \m{\lambda_2} pour lesquelles il existe des morphismes
semi-stables non stables, ainsi que des morphismes stables, sont les
$$\alpha_k = \frac{k}{n_1}, \ \ \ 1\leq k\leq \ {\rm Inf}(n,n_1-1).$$
Posons \ \m{m=1+{\rm Inf}(n,n_1-1)} et
$$\alpha_0=0, \ \ \ \alpha_m= \ {\rm Inf}(1,\frac{n+1}{n_1}).$$
Si \ \m{1\leq k\leq \ {\rm Inf}(n,n_1-1)}, et \ \m{\alpha_{k-1}<\lambda_2
<\alpha_k}, le morphisme \m{(f_1,f_2)} est stable si et seulement si 
\begin{itemize}
\item[--] 1 - Si l'image de $f_2$ est contenue dans
\ $\ko\ot N'_1\subset\ko\ot N_1$, on a \ $\dim(N'_1)\geq k$.
\item[--] Pour tout \ $(f'_1,f'_2)\in H.(f_1,f_2)$, si l'image de $f'_1$ est
contenue dans \ $\ko\ot N'_1\subset\ko\ot N_1$, on a \
$\dim(N'_1)\geq n_1-k+1$. 
\end{itemize}
 
\vskip 0.8cm 

\subsubsubsection{5.4.3.2}{Les constantes}

\medskip

Pour appliquer le th\'eor\`eme 5.6 on a besoin des constantes \m{c_0(1)},
\m{c'_0(1)}, \m{c_1(1)} et \m{c_2(1)}. On calcule ais\'ement que
$$c_0(1)=0, \ \ c'_0(1)=c_1(1)=\frac{n+1}{2}, \ \ c_2(1)=\frac{2n}{n^2+n-2}.$$

\vskip 0.8cm

\subsubsubsection{5.4.3.3}{Application du th\'eor\`eme 5.6, 1-}

\medskip

On obtient un quotient si
$$\frac{\lambda_2}{\lambda_1} \ > \ n+1,$$
c'est-\`a-dire si
$$\lambda_2 \ > \ \frac{n+1}{n+2}.$$
On obtient la seule vari\'et\'e de modules \m{M_{n_1}} si \ \m{n_1\leq n+1},
et aucune si \ \m{n_1>n+1}. 

\vskip 0.8cm

\subsubsubsection{5.4.3.4}{Application du th\'eor\`eme 5.6, 2-}

\medskip

On doit avoir
$$\lambda_2 \ > \ \frac{1}{n+2}+\frac{n(n+1)}{2n_1(n+2)} \ \ \ \ {\rm si } \
n_1\leq n+1,$$
$$\lambda_2 \ > \ 1-\frac{n+1}{2n_1} \ \ \ \ \ {\rm si } \ n_1>n+1.$$

Dans le premier cas, on obtient la construction de \m{M_k} si
$$k \ > \frac{n_1}{n+2}+\frac{n(n+1)}{2(n+2)}.$$
Ceci donne des vari\'et\'es de modules suppl\'ementaires si
$$n_1 \ \geq \ \lbrack\frac{n}{2}\rbrack + 2.$$

Dans le second cas, puisque \ \m{\lambda_2<\frac{n+1}{n_1}}, on doit avoir
$$n_1<\frac{3}{2}(n+1).$$
On obtient alors les vari\'et\'es de modules \m{M_k} pour
$$n_1-\frac{n+1}{2} \ < \ k \ \leq \ n+1.$$

\vskip 0.8cm

\subsubsubsection{5.4.3.5}{Application du th\'eor\`eme 5.6, 3-}
 
\medskip

On obtient les m\^emes vari\'et\'es de modules que pr\'ec\'edemment si
\ \m{n_1\leq n+1}. 

Si \ \m{n_1>n+1}, on sait construire le quotient si 
$$\lambda_2 \ > \ \frac{1}{2}+\frac{n}{2n^2-4}.$$
On sait dans ce cas construire \m{M_k} pour
$$n_1(\frac{1}{2}+\frac{n}{2n^2-4}) \ < \ k \ \leq n+1.$$
On obtient donc d'autres vari\'et\'es de modules de morphismes 
si \ \m{n+3\leq n_1<2n} \ et si $n$ est assez grand.

\vskip 0.8cm

\subsubsubsection{5.4.3.6}{Polarisations pathologiques. Cas o\`u 
\ \m{\lambda_2<1/2}}

\medskip

On suppose que \ \m{n_1} est pair : \ \m{n_1=2p}, et que \ \m{n_1\leq 2n+2}.
Soit \m{(z_1,\ldots,z_{n+1})} une base de \m{H^0(\ko(1))}. On consid\`ere le
morphisme
$$(f_1,f_2) : \ko(-2)\oplus\ko(-1)\lra\ko\ot\C^{2p}$$ 
o\`u \m{f_2}, \m{f_1} sont d\'efinis respectivement par les matrices
$$\left(\begin{array}{c}
z_1\cr .\cr .\cr .\cr z_p\cr 0\cr 0\cr .\cr .\cr .\cr 0
\end{array}\right), \ \ \ \left(\begin{array}{c}
z_2^2\cr .\cr .\cr .\cr z_{p+1}^2\cr z_1^2\cr z_1z_2\cr .\cr .\cr .\cr
z_1z_p\end{array}\right).$$ 
Alors, si \ \m{\lambda_2<1/2}, il est ais\'e de voir que \m{(f_1,f_2)} est
stable. Cependant son stabilisateur dans $G$ n'est pas r\'eduit \`a \m{\C^{*}}.
Il ne peut donc pas y avoir de quotient g\'eom\'etrique de \m{W^{s}_C} par
$G$ dans ce cas.

\newpage

\subsubsubsection{5.4.3.7}{Polarisations pathologiques. Cas o\`u 
\ \m{\lambda_2=1/2}}

\medskip

On suppose que $n$ est
pair et \ \m{n_1=n+2}. Soient \m{K_1}, \m{K_2} des sous-espaces vectoriels de
\m{H^0(\ko(1))} de dimension \m{\frac{n+2}{2}}, et \m{D} une droite de
\m{H^0(\ko(1))}. Soient
$$(z_1,\ldots,z_{\frac{n+2}{2}}), \ \ (z'_1,\ldots,z'_{\frac{n+2}{2}}), \ \ z$$
des bases de \m{K_1}, \m{K_2} et $D$ respectivement. Alors l'\'el\'ement
de \m{W_C} d\'efini par les matrices
$$\left(\begin{array}{c}
z_1\cr .\cr .\cr .\cr z_{\frac{n+2}{2}}\cr 0\cr .\cr .\cr .\cr
0\cr\end{array}\right)
 \ \ , \ \ \left(\begin{array}{c}
 0\cr .\cr .\cr .\cr 0\cr zz'_1\cr .\cr .\cr .\cr
zz'_{\frac{n+2}{2}}\cr\end{array}\right)$$
est semi-stable et sa $G$-orbite est ferm\'ee et ne 
d\'epend que de \m{K_1}, \m{K_2} et $D$. On la note \m{\phi(K_1,K_2,D)}. 
Remarquons que si \ \m{(K'_1,K'_2,D')\not = (K_1,K_2,D)}, on a
$$\phi(K'_1,K'_2,D') \ \not = \ \phi(K_1,K_2,D).$$
Une mutation de \m{\phi(K_1,K_2,D)} dans \m{W'_C} est un morphisme
$$\psi : \ko(-1)\ot\C^{n+2}\lra Q(-1)\oplus(\ko\ot\C^{n+2}).$$
L'adh\'erence de sa \m{G'}-orbite contient le morphisme somme directe des
morphismes
$$\psi_1 : \ko(-1)\lra\ko\ot\C^{n+2}, \ \ \psi_2 : \ko(-1)\lra\ko\ot\C^{n+2},
\ \psi_3 : \ko(-1)\ot\C^{n}\lra Q(-1)$$
d\'efinis respectivement par \m{K_1}, \m{K_2}, \m{D}. La $G'$-orbite de
ce morphisme est not\'ee \m{\psi(K_1,K_2,D)}. Notons qu'elle est contenue 
dans le compl\'ementaire dans \m{{W'}^{ss}_C} de l'ouvert contitu\'e
des orbites des morphismes mutations de morphismes de \m{W^{ss}_C}. On
a
$$\psi(K_1,K_2,D) \ = \ \psi(K_2,K_1,D).$$
Ceci prouve qu'on ne peut pas obtenir par des mutations indirectes un quotient
s\'eparant les orbites de \m{\phi(K_1,K_2,D)} et \m{\phi(K_2,K_1,D)}.

\vskip 2cm

{\bf Notes : } Ce texte reproduit l'article

{\em Quotients alg\'ebriques par des groupes non r\'eductifs et vari\'et\'es de 
modules de complexes.} International Journal of Mathematics Vol. 9 No 7 (1998) 
, 769-819.

avec quelques am\'eliorations dans la r\'edaction et la bibliographie. Les 
d\'efinitions de quotients dans 2.1 ont aussi \'et\'e modifi\'ees.

\vskip 1.5cm

\end{document}